\documentclass[12pt]{article}

\usepackage{latexsym,amsmath,amsthm,amssymb,fullpage}

\usepackage{color}

\usepackage{latexsym,amsmath,amsthm,amssymb,fullpage}

\usepackage[utf8]{inputenc}

\usepackage{amssymb}
\usepackage{epstopdf}
\usepackage{amsmath}
\usepackage{amsthm}
\usepackage{amsfonts}
\usepackage{epstopdf}

\newtheorem{theorem}{Theorem}

\newtheorem{proposition}[theorem]{Proposition}
\theoremstyle{definition}

\theoremstyle{remark}

\newcommand\beq{\begin {equation}}
\newcommand\eeq{\end {equation}}
\newcommand\beqs{\begin {equation*}}
\newcommand\eeqs{\end {equation*}}

\newcommand\N{{\mathbb N}}

\newcommand\R{{\mathbb R}}

\newcommand\E{{\mathbb E}}

\newcommand{\mathbbm}[1]{\text{\usefont{U}{bbm}{m}{n}#1}}

\begin {document}

$$   $$

\vskip 4mm
\begin{center}
\textbf{\LARGE On optimal matching of Gaussian samples III\\}

\vskip 10mm
\textit{\large Michel Ledoux, Jie-Xiang Zhu}\\
\vskip 3mm
{\large University of Toulouse, France, and Fudan University, China}\\
\end{center}
\vskip 7mm

\vskip 4mm

\begin {abstract}
This article is a continuation of the papers \cite {L17,L18} in which the optimal matching
problem, and the related rates of convergence of empirical measures for Gaussian
samples are addressed. A further step in both the dimensional
and Kantorovich parameters is achieved here,
proving that, given $X_1, \ldots, X_n$ independent random variables with common
distribution the standard Gaussian measure
$\mu$ on $\R^d$, $d \geq 3$, and $\mu_n \, = \, \frac 1n \sum_{i=1}^n \delta_{X_i}$ the associated
empirical measure,
$$
    \E  \big( \mathrm {W}_p^p (\mu_n , \mu )\big )
   \, \approx \,  \frac {1}{n^{p/d}}
$$
for any $1\leq p < d$, where $\mathrm {W}_p$ is the $p$-th Kantorovich metric.
The proof relies on the pde and mass transportation approach developed by
L.~Ambrosio, F.~Stra and D.~Trevisan in a compact setting.
\end {abstract}

\vskip 4mm

\section{Introduction and main results} \label {sec.1}
\setcounter{equation}{0}

Given $p \geq 1$, the Kantorovich distance (cf. \cite {V09} e.g.) between
two probability measures $\nu$ and $\mu$ on the Borel sets of $\R^d$ with a finite
$p$-th moment is defined by
\beq \label {eq.kanto}
\mathrm {W}_p(\nu ,\mu ) \, = \, \inf \bigg ( \int_{\R^d \times \R^d} |x-y|^p d\pi (x,y) \bigg)^{1/p}
\eeq
where the infimum is taken over all couplings $\pi $
on $\R^d \times \R^d$ with respective marginals $\nu $ and $\mu $. Here
$|x-y|$ denotes the Euclidean distance between $x$ and $y$ in $\R^d$.

Denote by $X_1, \ldots, X_n$, $n \geq 1$, independent random variables in $\R^d$
with common distribution $\mu$ and let
$$
\mu_n \, = \, \frac 1n \sum_{i=1}^n \delta_{X_i}
$$
be the empirical measure on the sample $(X_1, \ldots, X_n)$.
The question investigated here is the order of decay in $n $ of the expectations
\beq \label {eq.wp}
\E \big ( \mathrm {W}_p^p (\mu_n , \mu ) \big)
\eeq
when the random variables $X_1, \ldots, X_n$ are independent with the same
standard Gaussian distribution $d\mu(x) = e^{-|x|^2/2} \frac {dx}{(2\pi)^{d/2}}$ on $\R^d$.
The optimal matching problem would consist in the study of
$\E  ( \mathrm {W}_p^p (\mu_n , \nu_n ))$ for two independent samples
$X_1, \ldots, X_n$ and $Y_1, \ldots, Y_n$ with
$ \mu_n = \frac 1n \sum_{i=1}^n \delta_{X_i}$ and
$ \nu_n = \frac 1n \sum_{i=1}^n \delta_{Y_i}$, which is easily compared to \eqref {eq.wp}.

As a continuation of \cite {L17}, we refer to this article for more background and motivation
in the study of the optimal matching problem and rates of convergence in Kantorovich metrics
of empirical measures of Gaussian samples.
To introduce to the results of this work, we nevertheless recall the picture for the uniform
distribution on $[0,1]^d$, as well as the known results so far in the Gaussian (and more general) setting.

Throughout the paper $A \lesssim B$ between two real positive numbers $A$ and $B$
means that $A \leq C B$ where $C>0$ is either numerical or depends on $p$, $d$, but not
on anything else. In the same way, a sentence like ``$A$ is bounded from above by $B$" has the
same meaning. The equivalence sign $A \approx B$ indicates that $A \lesssim B$ and
$B \lesssim A$. In particular, these extended inequalities will then hold uniformly over $n \geq 1$.
Actually, $n$ may always be assumed to be larger than some fixed integer $n_0$,
large enough for obvious inequalities to hold true.

If $\mu$ is uniform on $[0,1]^d$, it holds true that (cf. e.g.~\cite {T14}),
for every $1 \leq p < \infty$,
\beq \label {eq.uniformp}
 \E \big ( \mathrm {W}_p^p (\mu_n, \mu) \big)  \, \approx \,
 \begin {cases}
  \frac {1}{n^{p/2}}   & \text{if \, $d=1$,} \\
    \big( \frac {\log n}{n}\big)^{p/2}    & \text{if \, $d=2$,} \\
   \frac {1}{n^{p/d}}  & \text {if \, $d\geq 3$.} \\
   \end {cases}
\eeq
The particular, and critical, case $d=2$ is the famous Ajtai-Koml\'os-Tusn\'ady theorem \cite {AKT84}.

Before addressing the Gaussian model, it is worthwhile mentioning that, by a simple
contraction argument (cf.~\cite {L17}), the expected cost
$ \E ( \mathrm {W}_p^p (\mu_n, \mu) ) $ when $\mu$ is uniform is bounded from
above by the corresponding quantity when $\mu$ is the standard Gaussian measure.
Hence, the rates in the uniform case provide lower bounds for the Gaussian model.
This comparison is used implicitly in the following descriptions.

Let now $\mu$ be the standard Gaussian measure on the Borel sets
of $\R^d$. It has been shown in \cite {BL16} that in dimension $d=1$,
\beq \label {eq.gaussian1}
 \E \big ( \mathrm {W}_p^p (\mu_n, \mu) \big)  \, \approx \,
 \begin {cases}
  \frac {1}{n^{p/2}}  & \text{if \, $1 \leq p < 2$,} \\
    \frac {\log \log n}{n}   & \text{if \, $p=2$,} \\
   \frac {1}{n (\log n)^{p/2}}   & \text {if \, $ p >2$.} \\
   \end {cases}
\eeq
With respect to \eqref {eq.uniformp}, it therefore appears that,
already in dimension one, the rates for ${p\geq 2}$ are rather sensitive
to the underlying distribution. The proof of \eqref {eq.gaussian1} in \cite {BL16} relies on monotone
rearrangement transport and explicit one-dimensional distributional inequalities.

In dimension $d=2$,
\beq \label {eq.gaussian2}
 \E \big ( \mathrm {W}_p^p (\mu_n, \mu) \big)  \, \approx \,
 \begin {cases}
  \Big (\frac {\log n}{n} \Big)^{p/2} & \text{if \, $1 \leq p < 2$,} \\
    \, \, \frac {(\log n)^2}{n}   & \text{if \, $p=2$.} \\
   \end {cases}
\eeq
Again, a specific new feature appears as $p=2$. As we will see, the case $p>2$ is essentially open.
The proof of the case $1 \leq p < 2$, and the upper bound for $p=2$, given in \cite {L17}
is based on the pde and mass transportation approach developed next, while the lower bound
for $p=2$ in \cite {T18} relies on the generic chaining ideas of \cite {T14} together with
a scaling argument. An alternate proof of this lower bound
using the pde-transportation method has been provided soon after in \cite {L18}.
When $p=1$, the upper bound
$ \E  ( \mathrm {W}_1 (\mu_n, \mu) )  \, \lesssim \, \sqrt {\frac {\log n}{n} }$
has been shown in \cite {Y92} to hold for distributions $\mu$ with a mild
moment assumption.

In higher dimension $d \geq 3$, a general bound using dyadic decompositions in the
spirit of the Ajtai-Koml\'os-Tusn\'ady theorem, and actually holding for
distributions $\mu$ with enough moments, has been obtained in \cite {DSS13,FG15} expressing that
\beq \label {eq.gaussianp<d2}
 \E \big ( \mathrm {W}_p^p (\mu_n, \mu) \big)  \, \lesssim \, \frac {1}{n^{p/d}}
\eeq
whenever $1\leq p < \frac d2$. For Gaussian samples, it is extended up to $p<2$ in dimension $3$
in \cite {L17}. It is also known from the works \cite {DSS13,FG15} that
\beq \label {eq.gaussianp>d2}
 \E \big ( \mathrm {W}_p^p (\mu_n, \mu) \big)
 		\,  \lesssim \,  \frac {1}{\sqrt n}
\eeq
when $p > \frac d2$, with an extra $\log n$ at the numerator when $p = \frac d2$.
The used methodology however will never produce anything better than this $\frac {1}{\sqrt n}$ rate,
and the upper bound \eqref {eq.gaussianp>d2} is actually far from the potential lower bound
deduced from \eqref {eq.uniformp}, and already not satisfactory when $d=1$.

The purpose of the work is to propose some progress in the understanding of the rates in this
Gaussian setting when $d \geq 3$ with the following statements.

\begin{theorem} \label {thm.main}
Let $X_1, \ldots, X_n$ be independent with common law the standard Gaussian distribution
$\mu$ on $\R^d$, $d \geq 3$,
and set $\mu_n = \frac{1}{n} \sum_{i=1}^{n} \delta_{X_i}$, $n \geq 1$. Then,
for $1 \leq p < d$,
$$
\E \big ( W_p^p( \mu_n, \mu ) \big) \, \approx \, \frac{1}{n^{p/d}} \, .
$$
\end{theorem}

In this range $1 \leq p < d$, $d \geq 3$, the rates for the Gaussian are therefore the same
as the ones for the compact uniform model.
The result extends \eqref {eq.gaussianp<d2} from $p < \frac d2$ to $p < d$.
This might look as only a small step, but it overcomes the
$\frac {1}{\sqrt n}$ rate and, as the proof will amply demonstrate, the amount
of work to reach this conclusion is rather significant. Due to the results
in \cite {DSS13,FG15,L17} when $1 \leq p < 2$ mentioned above,
only the values $2 \leq p < d$ have to be considered.

As identified by \eqref {eq.gaussian2} when $p=d=2$, the case $p=d$ might be of special interest.
We have been able to reach the following conclusion, not definitive however.

\begin{theorem} \label {thm.pd}
Let $X_1, \ldots, X_n$ be independent with common law the standard Gaussian distribution
$\mu$ on $\R^d$, $ d\geq 2$,
and set $\mu_n = \frac{1}{n} \sum_{i=1}^{n} \delta_{X_i}$, $n \geq 1$. Then,
$$
\E \big ( \mathrm {W}_d^d (\mu_n, \mu) \big)  \, \lesssim \, \frac {(\log n)^\kappa}{n}
$$
where
\beqs
\kappa \, = \,
 \begin {cases}
\frac{d^2 + 6d }{8} & \text{if \, $d=2,3$,} \\
   \frac {d^2}{2} - \frac {3d}{4}  & \text{if \, $d\geq 4$.} \\
   \end {cases}
\eeqs
\end{theorem}

The case $d=2$ recovers the claim in \eqref {eq.gaussian2} (although the proof
developed here to study every $d \geq 2$ is more involved than the ones in
\cite {L17} and \cite {T18}). The lower bound provided by the uniform model is $\frac 1n$ when
$d \geq 3$. A possible conjecture might be that
$$
\E  \big ( \mathrm {W}_d^d (\mu_n, \mu) \big) \, \approx \, \frac {(\log n)^{\frac d2}}{n}
$$
for $d \geq 3$. This is what is suggested as a lower bound in \cite {T18}.

At this point, we do not have any reasonable conjecture for $p >d$ ($\geq 2$).

Our proof of the preceding main results relies on the pde and transportation approach
developed by L.~Ambrosio, F.~Stra and D.~Trevisan \cite {AST19} towards the exact limit
for the two-dimensional uniform model in $\mathrm {W}_2$,
already exploited in the Gaussian case in \cite {L17}.
This approach relies on a heat kernel regularization argument together with transportation
bounds in terms of dual Sobolev norms. With respect to the compact case,
the Gaussian model involves the so-called Mehler kernel as underlying dynamics, which is
unbounded. One of the features of the work \cite {L17} was the introduction of a localization step in order
to take into account the infinite support of the Gaussian measure and the
associated Gaussian tails. A further main step is achieved
here by a randomization of the regularization time, together with estimates
on the Mehler kernel.

The note \cite {T18} by M. Talagrand gave a proof of the exceptional case $p=d=2$
in \eqref {eq.gaussian2} by means of a scaling argument (for more general distributions than the Gaussian).
It is certainly possible that the tools developed therein could lead to the
main results presented here, as well as answer some of the left open rates, in particular
in the case $p=d$. However, this note is rather difficult to grasp and we could not extract
further conclusions at this stage.
On the other hand, the pde and transportation method may be used to recover the main result
of \cite {T18} for $p=d=2$ as was shown in \cite {L18}
soon after \cite {T18}. For completeness and convenience, this proof is briefly recalled
here as Section~\ref {sec.7}.

Turning to the content of the paper,
Section~\ref {sec.2} collects several formulas and bounds on the Mehler kernel
as substitutes of uniform bounds in the compact case.
The subsequent paragraph briefly recalls, in this Gaussian context,
the transportation inequalities needed for the proofs already emphasized in \cite {AST19,L17}.
The proof of the main Theorem~\ref {thm.main} is divided in the two next sections, addressing
respectively the case $p=2$ and $p>2$. It could have been possible to present the proof
directly for $p\geq 2$, but the case $p=2$ is easier to handle and provides a good warm-up
for the general case. In Section~\ref {sec.6}, we discuss how the parameters may be adjusted
to reach Theorem~\ref {thm.pd}. As announced, the
last Section~\ref {sec.7} presents the pde and transportation proof
of the lower bound in \eqref {eq.gaussian2} first achieved by M.~Talagrand \cite {T18}
using combinatorial and scaling properties.

\section{Properties of the Mehler kernel} \label {sec.2}
\setcounter{equation}{0}

This short section collects basic and classical properties of the Mehler kernel and the associated
Ornstein-Uhlenbeck semigroup which will be of use throughout this work. The reference
\cite {BGL14} (and the bibliography therein) covers most of the claims emphasized below. With respect to
the compact case, the main issue here is that the Mehler kernel is
spacially unbounded, and so various
pointwise and integral controls have to be identified.

Let $d\mu(x) = e^{-|x|^2/2} \frac {dx}{(2\pi)^{d/2}}$ be the standard Gaussian measure
on the Borel sets of $\R^d$.
The Mehler kernel is given, for $t >0$, $x, y \in \R^d$, by
\beq \label {eq.mehler}
p_t(x,y) \, = \, p_t(y,x) \, = \,  \frac {1}{(1 - e^{-2t})^{d/2}}
		\exp \bigg ( - \frac {e^{-2t}}{2(1-e^{-2t})} \big [ |x|^2 + |y|^2 - 2e^t \, x \cdot y \big] \bigg).
\eeq
Here $|x|$ denotes the Euclidean length of $x \in \R^d$.
It holds true that $\int_{\R^d} p_t(x,y) d\mu(y) = 1$ for all $t>0$ and $x \in \R^d$.
The Mehler kernel satisfies besides the basic semigroup property with respect to $\mu$,
\beq \label {eq.mehlersemigroup}
\int_{\R^d} p_s(x,z) \, p_t(z,y) d\mu (z) \, = \, p_{s+t}(x,y)
\eeq
for all $s,t >0$ and $x,y \in \R^d$.

To ease the notation, set below $a = e^{-t} \in (0,1)$ in some of the statements.

When $x=y$,
\beq \label {eq.mehlerdiagonal}
p_t(x,x) \, = \, \frac {1}{(1 - a^2)^{d/2}} \, e^{\frac {a}{1+a} |x|^2}
		\, \leq \, \frac {1}{(1 - a^2)^{d/2}} \, e^{\frac { |x|^2}{2} } .
\eeq
It is easily seen besides that for all $x,y \in \R^d$,
\beq \label {eq.mehleruniform}
p_t(x,y)  \, \leq \, \frac {1}{(1 - a^2)^{d/2}} \, e^{\frac { |x|^2}{2} } .
\eeq
Combining \eqref {eq.mehleruniform}, \eqref {eq.mehlersemigroup} and \eqref {eq.mehlerdiagonal}
also shows that for every $q \geq 2$ and all $x \in \R^d$,
\beq \label {eq.mehlerintegral}
 \int_{\R^d} p_t(x,y)^q d\mu(y)
 	\, = \,  \int_{\R^d} p_t(x,y)^{q-2} p_t(x,y)^2 d\mu(y)
		\, \leq \, \frac {1}{(1 - a^2)^{(q-1)d/2}} \, e^{(q -1) |x|^2/2 } .
\eeq

The Mehler kernel generates the Ornstein-Uhlenbeck semigroup
\beq \label {eq.ou}
P_t f(x) \, = \, \int_{\R^d} f(y) \, p_t(x,y) d\mu (y)
		\, = \, \int_{\R^d} f \big ( e^{-t} x + \sqrt {1 - e^{-2t}} \, y \big ) d\mu (y)
\eeq
for all $t >0$, $ x \in \R^d$, and any suitable measurable function $f : \R^d \to \R$,
with the natural extension $P_0 = \mathrm {Id}$. The family ${(P_t)}_{t \geq 0}$ defines
a Markov semigroup, symmetric in $\mathrm {L}^2(\mu)$,
with infinitesimal generator $\mathrm {L} = \Delta - x \cdot \nabla$ for
which the integration by parts formula
\beq \label {eq.ibp}
\int_{\R^d} f (-\mathrm {L} g) d\mu \, = \, \int_{\R^d} \nabla f \cdot \nabla g \, d\mu
\eeq
holds true for every smooth functions $f,g : \R^d \to \R$. The spectrum of the operator
$\mathrm {L}$ is $\N$.

The semigroup ${(P_t)}_{t \geq 0}$ is a contraction in all $\mathrm {L}^p(\mu)$-spaces
with norms ${\| \cdot \|}_p$, $1 \leq p \leq \infty$.
The spectral gap induces an exponential decay for mean zero functions $f$ in $\mathrm {L}^2(\mu)$,
\beq \label {eq.exponentialdecay}
{\| P_t f \|}_2 \, \leq \, e^{-t} \,{\| f \|}_2, \quad t \geq 0.
\eeq
The hypercontractivity property on the other hand expresses that whenever $1 < p < q < \infty$
and $e^{2t} \geq \frac {q-1}{p-1}$,
\beq \label {eq.hypercontractivity}
{\| P_t f \|}_q \, \leq \, {\| f \|}_p.
\eeq
It will be convenient later on to combine the preceding smoothing properties in the following form:
for any $f : \R^d \to \R$ in $\mathrm {L}^p(\mu)$, $p \geq 2$, with mean zero,
\beq \label {eq.hypercontractivityp}
{\| P_tf\|}_p \, \leq \, C e^{-t/2} \, {\|  f \|}_p, \quad t \geq 0,
\eeq
where $C>0$ only depends on $p$. The decay $e^{-t/2}$ is far from optimal
but good enough for our purpose. For the proof, let $t_0>0$ be such that
$e^{t_0} = p-1 > 1$, so that by \eqref {eq.hypercontractivity}, for every $t \geq t_0$,
$$
{\| P_t f \|}_p  \, = \, {\big \| P_{t/2} (P_{t/2} f) \big \|}_p \, \leq \, {\| P_{t/2} f \|}_2.
$$
By the exponential decay in $\mathrm {L}^2(\mu)$-norm \eqref {eq.exponentialdecay},
${\| P_{t/2} f \|}_2 \leq e^{-t/2} {\| f\|}_2  \leq e^{-t/2} {\| f \|}_p$ and hence
$$
{\| P_t f \|}_p \, \leq \, e^{-t/2} \, {\|  f \|}_p.
$$
If $t \leq t_0$, ${\| P_t f \|}_p \leq {\| f \|}_p \leq C e^{-t/2} {\| f \|}_p$
with $C = e^{t_0/2}$. The claim \eqref {eq.hypercontractivityp} is established.

The following pseudo-Poincaré inequality is another useful tool:
for any $1 \leq p < \infty$, there exists
$C = C(p,d) >0$ such that for every smooth function $f : \R^d \to \R$ and every $0 < t \leq 1$,
\beq \label {eq.pseudopoincarep}
\int_{\R^d} |P_t f - f |^p d\mu  \, \leq \, C \, t^{p/2} \int_{\R^d} |\nabla f |^p d\mu .
\eeq
In case $p=2$, this inequality may be easily deduced spectrally.
For a proof in the general case, let $g: \R^d \to \R$ be smooth with ${\|g\|}_q \leq 1$ such that
$$
{\| P_t f - f \|}_p \, = \, \int_{\R^d} g (P_t f - f) d\mu
$$
where $\frac 1p + \frac 1q = 1$. By symmetry of the semigroup
${(P_t)}_{t \geq 0}$ and integration by parts
$$
\int_{\R^d} g (P_t f - f) d\mu
	\, = \, \int_0^t \int_{\R^d} g \, \mathrm {L}P_s f \, d\mu \, ds
	\, = \, - \int_0^t \int_{\R^d} \nabla P_sg \cdot \nabla f \, d\mu \, ds.
$$
The integral representation \eqref {eq.ou} and one more integration by parts indicate that,
at any point $x \in \R^d$,
\beqs \begin {split}
 \nabla P_sg (x)
 	&\, = \, e^{-s} \int_{\R^d} \nabla g \big ( e^{-s} x + \sqrt {1 - e^{-2s}} \, y \big ) d\mu (y) \\
	&\, = \, \frac {e^{-s}}{\sqrt {1-e^{-2s}}}
		 \int_{\R^d} y \, g \big ( e^{-s} x + \sqrt {1 - e^{-2s}} \, y \big ) d\mu (y). \\
\end {split} \eeqs
Hence, by Hölder's inequality,
\beqs \begin {split}
\int_{\R^d} \nabla P_sg \cdot \nabla f \, d\mu
	&\, = \, \frac {e^{-s}}{\sqrt {1-e^{-2s}}}
		 \int_{\R^d} \! \int_{\R^d} y \cdot \nabla f(x) \, g \big ( e^{-s} x + \sqrt {1 - e^{-2s}} \, y \big )
		 d\mu (x) d\mu (y) \\
	&\, \leq \, \frac {e^{-s}}{\sqrt {1-e^{-2s}}}
		 \bigg (\int_{\R^d} \! \int_{\R^d}\big | y \cdot \nabla f(x) \big|^p d\mu (x) d\mu (y)
		 \bigg)^{1/p} \\
\end {split} \eeqs
where it is used that
$$
\int_{\R^d} \! \int_{\R^d} \big | g \big ( e^{-s} x + \sqrt {1 - e^{-2s}} \, y \big) \big |^q d\mu(x) d\mu(y)
		\, = \, \int_{\R^d} |g|^q d\mu \,  \leq \, 1.
$$
Finally, after partial integration in $d\mu(y)$,
$$
\int_{\R^d} \int_{\R^d} \big |y \cdot \nabla f(x) \big |^p  d\mu (x) d\mu (y)
 	\, = \,  C \int_{\R^d}  |\nabla f |^p  d\mu
$$
where $C>0 $ only depends on $p$. The claim \eqref {eq.pseudopoincarep} then easily follows.

We will make use of an important and delicate property, the Riesz transform bounds.
In this Gaussian setting, there were established by P. A.~Meyer \cite {M84}
(see also \cite {B87}) and express that, for every $1 < p < \infty$
and every $f : \R^d \to \R$ in the suitable domain,
\beq \label {eq.riesz}
\int_{\R^d} |\nabla f|^p d\mu \, \approx \, \int_{\R^d} \big | (- \mathrm {L})^{-1/2} f \big|^p d\mu
\eeq
($\approx $ only depending on $p$)
where $(- \mathrm {L})^{-1/2}$ is defined spectrally on mean zero functions $f$,
for example by the classical formula
\beq \label {eq.spectral}
(- \mathrm {L})^{-1/2} \, = \, \frac {1}{\sqrt \pi} \int_0^\infty \! \frac {1}{\sqrt s} \, P_s  \, ds.
\eeq

Finally, some technical tools related to energy estimates will be requested.
The reverse Poincaré inequality for the Gaussian measure
(\cite [Theorem 4.7.2] {BGL14}) expresses that for every
Borel set $A$ in $\R^d$ and every $s>0$,
$$
|\nabla P_s\mathbbm {1}_A|^2
 	\, \leq \, \frac{1}{e^{2s} - 1} \, \big [P_s(\mathbbm {1}_A^2) - (P_s \mathbbm {1}_A)^2\big]
 	\, \leq \,  \frac{1}{e^{2s} - 1} \, .
$$
Combining with \eqref {eq.pseudopoincarep} for $f=P_s\mathbbm {1}_A$, for every $s >0$,
$t \in (0,1)$ and $ p \geq 1$,
\beq \begin {split} \label {eq.pseudopoincare1}
\int_{\R^d} \big |P_t (P_s\mathbbm {1}_A) - P_s \mathbbm {1}_A \big|^p d\mu
 & \, \leq \, \frac{C \, t^{p/2}}{(e^{2s} - 1)^{(p-1)/2}} \int_{\R^d} |\nabla P_s \mathbbm {1}_A | d\mu \\
 & \, \leq \, \frac{C \, t^{p/2}\, e^{-s}}{(e^{2s} - 1)^{(p-1)/2}} \, \mu (\partial A)\\
\end {split} \eeq
where, in the last step, it is used that
$\nabla P_s = e^{-s} P_s \nabla$ and
$\limsup_{\varepsilon \to 0} \int_{\R^d} |\nabla P_\varepsilon \mathbbm {1}_A| d\mu
 \leq \mu(\partial A)$ for a measurable subset $A$ of $\R^d$
 with smooth boundary $\partial A$.

\section{Mass transportation bounds} \label {sec.3}
\setcounter{equation}{0}

This short paragraph presents the key functional analytic tool to bound
Kantorovich distances in this context. It was already put forward in \cite {L17}
(see also \cite {P18,S15}).

\begin{proposition} \label{prop.sobolev}
Let $p \geq 1$. For any $d\nu = f d\mu$ with $f-1$ in $\mathrm {H}^{-1, p}(\mu)$, we have
$$
\mathrm {W}_p(\nu, \mu) \, \leq \,  p \, {\| f-1 \|}_{\mathrm {H}^{-1, p}(\mu)}
$$
where the $\mathrm {H}^{-1, p}(\mu)$ negative Sobolev norm is defined by
$$
{\| g \|}_{\mathrm {H}^{-1, p}(\mu)}
	\, = \, \bigg( \int_{\R^d} \big | \nabla ((- \mathrm {L})^{-1} g) \big |^p d\mu \bigg)^{1/p}
$$
for any $g : \R^d \to \R$ with mean zero such that the right-hand side makes sense.
\end{proposition}

It may be emphasized that in the particular case $p=2$,
\beq \label {eq.sobolev2}
{\| g \|}^2_{\mathrm {H}^{-1, 2}(\mu)}  \, = \,
  2 \int_{0}^{\infty} \! \int_{\R^d} (P_s g)^2 d\mu \, ds.
\eeq
Indeed, since $(-\mathrm {L})^{-1} = \int_0^\infty P_s ds$, by integration by parts,
$$
\int_{\R^d} \big | \nabla ((- \mathrm {L})^{-1} g) \big |^2 d\mu
 	\, = \, \int_{\R^d} g (- \mathrm {L})^{-1} g  d\mu
	\, = \, \int_0^\infty \! \int_{\R^d} g \, P_s g  \,d\mu \, ds
$$
from which the claim follows from the symmetry of $P_s$.

When $p\not= 2$, it will be necessary to rely on the Riesz transform bound \eqref {eq.riesz},
for $1 < p < \infty$,
after which Proposition~\ref {prop.sobolev} takes the form
\beq \label {eq.sobolevriesz}
\mathrm {W}_p^p(\nu, \mu) \, \lesssim \,
      \int_{\R^d} \big | (- \mathrm {L})^{-1/2} (f-1) \big |^p d\mu .
\eeq

\section{The case $2 = p < d$} \label {sec.4}
\setcounter{equation}{0}

We therefore address here the proof of Theorem~\ref {thm.main} for $2 = p < d$.
Although the proof of this case is actually contained in the more general section $2 \leq p < d$,
it is easier to access due to the semigroup representation \eqref {eq.sobolev2} of the negative Sobolev norm.

The first step in the investigation is the localization argument introduced in \cite {L17}.
For $R > 0$, let $d\mu^R = \frac{1}{\mu(B_R)} \mathbbm {1}_{B_R} d\mu$ where $B_R$ is the Euclidean ball centered at 0 with radius $R$. Define independent random variables $X^R_i$, $i = 1, \ldots, n$, with common distribution $\mu^R$ by
\begin{align*}
X^R_i \, = \,  \left\{ \begin{array}{cc}
X_i  & \textrm{if    } X_i \in B_R, \\
Z_i  & \textrm{if    } X_i \notin B_R,
\end{array}
\right.
\end{align*}
where $Z_1, \ldots, Z_n$ are independent with distribution $\mu^R$, independent of the $X_i$'s.
Setting $\mu^R_n = \frac{1}{n} \sum_{i=1}^{n} \delta_{X^R_i}$, by definition of the coupling,
$$
\mathrm {W}^2_2(\mu_n, \mu^R_n)
	\, \leq \, \frac{1}{n} \sum_{i=1}^{n} |X_i - X^R_i|^2
	\, \leq \, \frac{4}{n} \sum_{i=1}^{n} |X_i|^2 \mathbbm {1}_{\{|X_i| \geq R\}}.
$$
Therefore
$$
\E \big ( \mathrm {W}_2^2( \mu_n, \mu^R_n ) \big) \,  \leq \,  4 \int_{\{|x| > R \}} |x|^2 d\mu.
$$
Since $ \int_{\{|x| > R \}} |x|^2 d\mu = C_d \int_R^{\infty} r^{d+1} e^{-r^2/2} dr$
is of the order of $R^de^{-R^2/2}$ as $R \to \infty$,
choose $R = \sqrt{2c\log n}$ for some $c \in (\frac{2}{d}, 1)$ so that
\beq  \label{eq.localization2}
\E \big ( \mathrm {W}_2^2( \mu_n, \mu^R_n ) \big)  \, \lesssim \, \frac {1}{n^{2/d}} \, .
\eeq
As a result, the investigation is concentrated on the study of
$\E ( \mathrm {W}_2^2( \mu_n^R, \mu ) )$.
Note furthermore that $\mu(B_R) \geq \frac 12$ for $n$ large enough so that this
normalization factor may essentially be neglected throughout the investigation. For the further
developments, it will be convenient to refer to a random variable $X$ with distribution $\mu$
and to $X^R$ with distribution $\mu^R$.

A new step with respect to the former investigations is
the introduction of a randomized regularization time by means of the decomposition
of $B_R$ as the union of $m$ annuli
$$
D_k  \,= \, \big \{ x \in \R^d ; \, r_{k-1} \leq |x| < r_k \big \}, \quad k = 1, \ldots, m,
$$
where $0 = r_0 < r_1 < \cdots <r_m = R$ with $r_k = \sqrt {k}$,
$k=1, \ldots, m$. In particular ${m = R^2 = 2c \log n}$.
Define then a map $T : B_R \to (0, 1)$
as $T(x) = t_k$ if $x \in D_k$ where the $ 0 < t_1 < \cdots <  t_m < 1$ will be specified later.

For this map $T$, consider then the (random) probability density
$$
f(y) \, = \,  f^{R, T}_n(y) = \frac{1}{n} \sum_{i=1}^{n} p_{T(X^R_i)} (X^R_i, y), \quad y \in \R^d,
$$
and set $d\mu^{R, T}_n = f^{R, T}_n d\mu$. By convexity of the Kantorovich
metric $\mathrm {W}_2^2$ \cite [Theorem 4.8] {V09}
and the representation formula
for the Ornstein-Uhlenbeck semigroup \eqref {eq.ou}, conditionally on the $X_i^R$'s,
\beqs \begin {split}
\mathrm {W}_2^2( \mu^R_n, \mu^{R, T}_n )
& \, \leq \, \frac{1}{n} \sum_{i=1}^{n} \mathrm {W}_2^2 \big ( \delta_{X^R_i},
		p_{T(X^R_i)} (X^R_i, \cdot)d\mu \big)\\
& \, = \, \frac{1}{n} \sum_{i=1}^{n} \int_{\R^d} |X^R_i - y|^2 \, p_{T(X_i^R)} (X^R_i, y) d\mu(y)\\
& \, = \,  \frac 1n \sum_{i=1}^n \Big [ \big(1 - e^{-T(X_i^R)}\big)^2 |X^R_i|^2
		+ d \big (1 - e^{-2T(X_i^R)} \big) \Big] \\
&  \, \lesssim \,  \frac 1n \sum_{i=1}^n  \big [ T(X_i^R)^2 |X^R_i|^2 + T (X_i^R) \big]
\end {split} \eeqs
where we used in the last step that $ T(X^R_i) \in (0, 1)$.
Averaging over the $X_i^R$'s, and dropping $X^R$ in $T(X^R)$ for simplicity,
\beq \begin {split} \label {eq.regularization2}
\E \big (\mathrm {W}_2^2( \mu^R_n, \mu^{R, T}_n ) \big)
 & \, \lesssim \,  \E \big( T^2 |X^R|^2 \big) + \E (T) \\
  & \, \lesssim \, (t_m R^2+ 1) \,  \E(T)
  \, = \, ( t_m R^2 + 1)  \sum_{k=1}^{m} t_k \mu(D_k). \\
\end {split} \eeq

From the localization \eqref {eq.localization2} and regularization
\eqref {eq.regularization2} arguments, and the triangle inequality for
$\mathrm {W}_2$, we have therefore obtained at this stage that
\beq \label {eq.localizationregularization2}
\E \big (\mathrm {W}_2^2( \mu_n, \mu ) \big)
		\, \lesssim \,  \E \big (\mathrm {W}_2^2( \mu_n^{R,T}, \mu ) \big)
		+ \frac {1}{n^{2/d}} + ( t_m R^2 + 1)  \sum_{k=1}^{m} t_k \mu(D_k).	
\eeq

From here, we thus concentrate on
$\E(\mathrm {W}_2^2( \mu^{R, T}_n, \mu ))$ for which we make use of Proposition~\ref {prop.sobolev}
and \eqref {eq.sobolev2} with
$$
g \, = \, g(y) \, = \,   f^{R, T}_n(y) - 1
	\, = \, \frac{1}{n} \sum_{i=1}^{n} \big [p_{T(X_i^R)} (X^R_i, y) - 1\big],
		\quad y \in \R^d,
$$
to get that
$$
\E \big (\mathrm {W}_2^2( \mu^{R, T}_n, \mu ) \big)
  \, \leq \, 4 \, \E \big({ \| f^{R, T}_n - 1 \|}^2_{\mathrm {H}^{-1, 2}(\mu)} \big)
   \, = \, 8 \, \int_{0}^{\infty} \! \int_{\R^d} \E \big( (P_s g)^2 \big) d\mu \, ds.
$$
In order to develop probabilistic arguments, it is
convenient to center the elements $p_{T(X_i^R)} (X^R_i, y)$ in the definition of $g$.
Write therefore, for every $y \in \R^d$,
$$
g(y) \, = \, \frac 1n \sum_{i=1}^n \big [p_{T(X_i^R)} (X^R_i, y)
		- \E (p_{T(X_i^R)} (X^R_i, y) ) \big]
	+ \E \big (p_T(X^R, y) \big) - 1  \, = \,  \widetilde{g}(y) + \phi(y)
$$
where $\phi (y) = \E (p_T(X^R, y)) - 1$. Recall that here $T = T(X^R)$.
For every $s>0$, $ \E ( (P_s g)^2 ) \leq 2\, \E ((P_s \widetilde{g})^2 ) + 2(P_s \phi)^2$
so that
$$
\E \big (\mathrm {W}_2^2( \mu^{R, T}_n, \mu ) \big)
   \, \leq \, 16 \, \int_{0}^{\infty} \! \int_{\R^d} \E \big( (P_s \widetilde {g})^2 \big) d\mu \, ds
   	+ 16 \, \int_{0}^{\infty} \! \int_{\R^d} (P_s \phi)^2 d\mu \, ds .
$$
Now, since
$$
P_s \widetilde {g} (y) \, = \,  \frac 1n \sum_{i=1}^n \big [p_{s+T(X_i^R)} (X^R_i, y)
		- \E (p_{s+T(X_i^R)}(X^R_i, y) ) \big]
$$
and since the $X^R_i$'s, $i= 1, \ldots, n$, are independent and identically distributed,
for each $y \in \R^d$,
$$
\E \big( (P_s \widetilde {g})^2 \big)
	 \, = \,  \frac{1}{n} \, \E \Big ( \big [  p_{s+T}(X^R, y) - \E(p_{s+T}(X^R, y)) \big]^2  \Big) .
$$
But $ \E(p_{s+T}(X^R, y))  = P_s \phi (y)  + 1$, and as we aim to control
$ \int_{0}^{\infty} \! \int_{\R^d} (P_s \phi)^2 d\mu \, ds$, we may as well replace back
$ \E(p_{T+s}(X^R, y))$ by $1$ in the latter, to obtain that
\beq \label {eq.2termdecomposition2}
\E \big (\mathrm {W}_2^2( \mu^{R, T}_n, \mu ) \big)
   \, \lesssim \,  \, \frac 1n \int_{0}^{\infty} \! \int_{\R^d}
   \E \big ( \big [  p_{s+T}(X^R, y) - 1 \big]^2  \big) d\mu(y) ds
   	+   \int_{0}^{\infty} \! \int_{\R^d} (P_s \phi)^2 d\mu \, ds .
\eeq

The most important term on the right-hand side of \eqref {eq.2termdecomposition2} is the
first one on which we concentrate next.
Divide the integral in $s$ according as $s \in (0,1)$ or $s \in (1,\infty)$ so to get
\beq \begin {split} \label {eq.2termdecomposition21}
 \frac 1n  \int_{0}^{\infty} \!   \int_{\R^d}
   		\E   \big ( \big [   p_{s+T}& (X^R,  y) - 1  \big]^2  \big) d\mu(y) ds \\
		& \, = \, \frac 1n \int_{0}^1 \! \int_{\R^d}
   \E \big ( \big [  p_{s+T}(X^R, y) - 1 \big]^2  \big) d\mu(y) ds  \\
  & \quad \, \, \, + \frac 1n \int_1^{\infty} \! \int_{\R^d}
   \E \big ( \big [  p_{s+T}(X^R, y) - 1 \big]^2  \big) d\mu(y) ds. \\
\end {split} \eeq
Recall in addition that by definition of $T$,
$$
\E \big ( \big [  p_{s+T}(X^R, y) - 1 \big]^2  \big)
	 \, = \frac {1}{\mu(B_R)}
	\sum_{k=1}^{m} \int_{D_k} \big [ p_{s + t_k}(x, y) - 1 \big]^2 d\mu(x) .
$$
The first piece on the right-hand side of \eqref {eq.2termdecomposition21} may then
be simply upper-bounded as
$$
\frac 1n \int_{0}^1 \! \int_{\R^d}
   \E \big ( \big [  p_{s+T}(X^R, y) - 1 \big]^2  \big) d\mu(y) ds
   \, \lesssim \, \frac 1n  \sum_{k=1}^m \int_{D_k} \int_{0}^1 \! \int_{\R^d}
    p_{s+t_k}(x, y)^2  d\mu(y) ds d\mu(x).
$$
By \eqref {eq.mehlersemigroup}, \eqref {eq.mehlerdiagonal} and integration in $s$,
for every $x \in \R^d$ and $k = 1, \ldots, m$,
$$
\int_{0}^1 \! \int_{\R^d}   p_{s+t_k}(x, y)^2   d\mu(y)ds
	\, = \, \int_{0}^1  p_{2(s+t_k)}(x, x)  ds
		\, \lesssim \, \frac {1}{t_k^{(d/2)-1}} \, e^{|x|^2/2} .
$$
Hence, since $\lambda (D_k) \approx  k^{(d/2)-1}$,
 \beq \begin {split}
\label {eq.decomposition2firstterm}
\frac 1n \int_{0}^1 \! \int_{\R^d}
   \E \big ( \big [  p_{s+T}(X^R, y) - 1 \big]^2  \big) d\mu(y) ds
   & \, \lesssim \, \frac 1n  \sum_{k=1}^m \int_{D_k} \frac {1}{t_k^{(d/2)-1}} \, e^{|x|^2/2} d\mu(x) \\
  &  \, \lesssim \, \frac 1n  \sum_{k=1}^m \frac {\lambda (D_k)}{t_k^{(d/2)-1}} \\
  &  \, \lesssim \, \frac 1n  \sum_{k=1}^m  \Big ( \frac {k}{t_k} \Big)^{(d/2)-1} . \\
\end {split} \eeq
On the other hand, towards the second piece on the right-hand side of \eqref {eq.2termdecomposition21},
by the exponential decay \eqref {eq.exponentialdecay} in $\mathrm {L}^2(\mu)$,
\beq \begin {split} \label {eq.decomposition2secondterm}
\frac 1n \int_1^\infty \! \int_{\R^d}
   \E \big ( \big [   p_{s+T}(&X^R, y) -   1 \big]^2  \big) d\mu(y) ds  \\
  & \, = \, \frac 1n  \,\frac {1}{\mu(B_R)} \sum_{k=1}^m \int_{D_k} \int_1^\infty \! \int_{\R^d}
   \big [ p_{s+t_k}(x, y)-  1 \big]^2  d\mu(y) ds d\mu(x) \\
 & \, \lesssim \, \frac 1n  \sum_{k=1}^m \int_{D_k} \int_0^\infty \! e^{-2s} \int_{\R^d}
   \big [ p_1(x, y)- 1 \big]^2  d\mu(y) ds d\mu(x) \\
	 & \, \lesssim \, \frac 1n  \sum_{k=1}^m \lambda (D_k)
	 \, = \, \frac 1n  \, \lambda (B_R)
\end {split} \eeq
where it is used again that $ \int_{\R^d} p_1(x, y)^2  d\mu(y) = p_2 (x,x) \lesssim e^{|x|^2/2}$.

Summarizing \eqref {eq.decomposition2firstterm} and \eqref {eq.decomposition2secondterm}
in \eqref {eq.2termdecomposition21}, it follows that
$$
\frac 1n \int_{0}^{\infty} \! \int_{\R^d}
   \E \big ( \big [  p_{s+T}(X^R, y) - 1 \big]^2  \big) d\mu(y) ds
		\, \lesssim \, \frac 1n \sum_{k=1}^{m} \Big ( \frac {k}{t_k} \Big)^{(d/2)-1}
				+ \frac {1}{n^{2/d}}
$$
(recall that $d > 2$).
Collecting this estimate in \eqref {eq.2termdecomposition2}
and \eqref {eq.localizationregularization2}, it holds true that
\beq \begin {split} \label {eq.intermediate2}
\E \big (\mathrm {W}_2^2( \mu_n, \mu ) \big)
	& \, \lesssim \,
	\frac 1n \sum_{k=1}^{m} \Big ( \frac {k}{t_k} \Big)^{(d/2)-1}
	+ (t_mR^2 + 1) \sum_{k=1}^m t_k \mu(D_k) \\
	& \quad \, \,  \, +  \frac {1}{n^{2/d}} + \int_{0}^{\infty} \! \int_{\R^d} (P_s \phi)^2 d\mu \, ds. \\
\end {split} \eeq

The choice of the $t_k$'s is now determined by optimization between the first two terms on the
right-hand side of \eqref {eq.intermediate2}. Using that
$\mu(D_k) \lesssim k^{(d/2)-1} e^{-k/2}$, set
$t_k = n^{-2/d} e^{k/d}$, $k= 1, \ldots, m$. In particular,
$$
t_k \, \leq \, t_m  \, = \, \frac {1}{n^{2(1-c)/d}}
$$
for every $k=1,\ldots, m = R^2 = 2c\log n$, $c \in (\frac 2d,1)$.
For these values, and since $d>2$, it follows that
$$
\frac 1n \sum_{k=1}^{m}  \Big ( \frac {k}{t_k} \Big)^{(d/2)-1}
	+ (t_m R^2 + 1) \sum_{k=1}^m t_k \mu(D_k) \, \lesssim \, \frac {1}{n^{2/d}} \, .
$$
Therefore \eqref {eq.intermediate2} yields
\beq \label {eq.intermediate21}
\E \big (\mathrm {W}_2^2( \mu_n, \mu ) \big)
	\, \lesssim \,  \frac {1}{n^{2/d}} + \int_{0}^{\infty} \! \int_{\R^d} (P_s \phi)^2 d\mu \, ds .
\eeq

We are left with the study of the centering term
$\int_{0}^{\infty} \! \int_{\R^d} (P_s \phi)^2 d\mu \, ds$. To this task, recall that
\beqs \begin {split}
P_s \phi
& \, = \, \E \big ( p_{s + T} (X^R, \cdot) - 1 \big ) \\
& \, = \,  \frac{1}{\mu(B_R)} \sum_{k = 1}^{m} \big [P_{s+ t_k}\mathbbm {1}_{D_k} - \mu(D_k)\big]\\
& \, = \,  \frac{1}{\mu(B_R)} \sum_{k = 1}^{m} \big [P_{s+t_k}\mathbbm {1}_{D_k}
		- P_s \mathbbm {1}_{D_k} \big ]
			+ \frac{1}{\mu(B_R)} \sum_{k = 1}^{m} \big [ P_s \mathbbm {1}_{D_k} - \mu(D_k) \big ]\\
& \, = \, \frac{1}{\mu(B_R)} \sum_{k = 1}^{m}
	\big [P_{s+ t_k}\mathbbm {1}_{D_k} - P_s \mathbbm {1}_{D_k} \big]
	 + \frac{1}{\mu(B_R)} \big [ P_s \mathbbm {1}_{B_R} - \mu(B_R) \big ].
\end{split} \eeqs
Hence, for every $s>0$,
$$
\int_{\R^d} (P_s \phi)^2 d\mu
	\, \lesssim \, \int_{\R^d} \bigg [ \sum_{k = 1}^{m}
	\big [P_{s+t_k}\mathbbm {1}_{D_k} - P_s \mathbbm {1}_{D_k} \big ] \bigg ]^2 d\mu
	 + \int_{\R^d} \big [  P_s \mathbbm {1}_{B_R} - \mu(B_R)  \big ]^2 d\mu
$$
and we treat separately the two expressions on the right-hand side. The second one is easy since
by the exponential decay \eqref {eq.exponentialdecay} in $\mathrm {L}^2(\mu)$,
$$
\int_{\R^d} \big [  P_s \mathbbm {1}_{B_R} - \mu(B_R)  \big ]^2 d\mu
	\, \leq \, e^{-2s} \big [ 1 - \mu (B_R) \big] \, \lesssim \, \frac {e^{-2s}}{n^{2/d}}
$$
by the choice of $R$.
For the first one, by the triangle inequality and \eqref {eq.pseudopoincare1},
\beqs \begin {split}
\int_{\R^d} \bigg [ \sum_{k = 1}^{m}
	\big [P_{s+t_k}\mathbbm {1}_{D_k} - P_s \mathbbm {1}_{D_k} \big ] \bigg ]^2 d\mu
	& \, \lesssim \, \bigg ( \sum_{k=1}^{m} \Big[ \frac{t_k \, e^{-s}}{\sqrt{e^{2s} - 1}}
			\, \mu(\partial D_k) \Big]^{1/2} \bigg )^2 \\
	& \, \lesssim \, \frac{e^{-s}}{\sqrt{e^{2s} - 1}}
	\bigg ( \sum_{k=1}^{m} \sqrt{t_k} \, \sqrt{ \mu(\partial D_k) } \bigg)^2.
\end{split} \eeqs
Combining the preceding,
$$
\int_{0}^{\infty} \! \int_{\R^d} (P_s \phi)^2 d\mu \, ds
	\, \lesssim \, \bigg ( \sum_{k=1}^{m} \sqrt{t_k} \, \sqrt{ \mu(\partial D_k) } \bigg)^2
		+ \frac {1}{n^{2/d}} \, .
$$	
Since $\mu( \partial D_k) \lesssim  k^{(d-1)/2} e^{-k/2}$,
the choice of $t_k = n^{-2/d} e^{k/d}$, $k = 1, \ldots, m$,
easily shows that this centering contribution
is at most $\frac {1}{n^{2/d}}$. Inserting this claim into \eqref {eq.intermediate21}
concludes the proof of the theorem for $2 = p < d$.

\section{The case $2 \leq p < d$} \label {sec.5}
\setcounter{equation}{0}

The pattern of the proof will be similar to the one of Section~\ref {sec.4}
but with significant increase of the technicalities
since \eqref {eq.sobolev2} is now more available and the arguments go through the more involved
Riesz transform bound \eqref {eq.sobolevriesz}. The scheme of proof is then similar to the
one developed in \cite {L17} in the compact case \eqref {eq.uniformp} but,
again, unboundedness of the Mehler kernel requires several delicate estimates.

As in the preceding section for $p=2$, we
truncate in the same way on a ball $B_R$ with $R = \sqrt{ 2c\log n}$ for some $c \in (\frac pd, 1)$
for which we get similarly that
\beq  \label{eq.localizationp}
\E \big ( \mathrm {W}_p^p( \mu_n, \mu^R_n ) \big)  \, \lesssim \, \frac {1}{n^{p/d}} \, .
\eeq
We decompose again $B_R$ as the union of $m$ annuli $D_k = \{x \in \R^d ; r_{k-1} \leq |x| < r_k \}$,
$k = 1, \ldots, m$, where $0 = r_0 < r_1 < \cdots <r_m = R$, and consider
also the map $T : B_R \to (0, 1)$ defined
by $T(x) = t_k$ if $x \in D_k$. We will use the same choices $r_k = \sqrt {k}$
and $t_k = n^{-2/d} e^{k/d}$, $k = 1, \ldots, m$. In particular, for the further purposes,
it is important to notice again that
\beq \label {eq.tm}
t_k \, \leq \,  t_m \, \leq \,  \frac {1}{n^{2(1-c)/d}} \, , \quad k= 1, \ldots, m,
\eeq
small enough therefore for a number of subsequent small issues (recall that $n \geq n_0$ is assumed large enough throughout the investigation).

For this map $T : B_R \to (0, 1)$, set similarly
$$
f(y) \, = \,  f^{R, T}_n(y) \, = \,  \frac{1}{n} \sum_{i=1}^{n} p_{T(X^R_i)} (X^R_i, y), \quad y \in \R^d,
$$
and $d\mu^{R, T}_n = f^{R, T}_n d\mu$. By convexity \cite [Theorem 4.8] {V09}
and the representation formula
for the Ornstein-Uhlenbeck semigroup \eqref {eq.ou}, conditionally on the $X_i^R$'s,
\beqs \begin {split}
\mathrm {W}_p^p( \mu^R_n, \mu^{R, T}_n )
& \, \leq \,  \frac{1}{n} \sum_{i=1}^{n}
		 \mathrm {W}_p^p \big ( \delta_{X^R_i}, p_{T(X_i^R)} (X^R_i, \cdot)d\mu \big )\\
&  \, = \,  \frac{1}{n} \sum_{i=1}^{n} \int_{\R^d} |X^R_i - y|^p \, p_{T(X_i^R)} (X^R_i, y) d\mu(y)\\
& \, \lesssim \,   \frac{1}{n} \sum_{i=1}^{n}
	\Big [ \big (1 - e^{-T(X_i^R)} \big)^p |X^R_i|^p
		+ \big (1 - e^{-2T(X_i^R) } \big)^{p/2} \Big ] .
\end {split} \eeqs
Noticing that $ T(X^R_i) \in (0, 1)$, it follows after taking expectation that
\beq  \label {eq.regularizationp}
\E \big ( \mathrm {W} _p^p( \mu^R_n, \mu^{R, T}_n ) \big )
 	  \, \lesssim \, \E \big (  T^p |X^R|^p \big) + \E ( T^{p/2})
	  \, \lesssim \, (t_m^{p/2} R^p + 1) \, \E ( T^{p/2})
 \eeq
where we recall that $ T = T(X^R)$.

From \eqref {eq.localizationp}, \eqref {eq.regularizationp} and the triangle inequality,
and the fact \eqref {eq.tm} that $t_m \, \leq \,  \frac {1}{n^{2(1-c)/d}}$, we
have therefore obtained at this point that
\beq  \begin {split} \label {eq.localizationregularizationp}
\E \big ( \mathrm {W} _p^p( \mu_n, \mu ) \big )
 	& \, \lesssim \,  \E \big ( \mathrm {W} _p^p( \mu_n^{R,T}, \mu ) \big )
 		+ \frac {1}{n^{p/d}} +  \sum_{k=1}^{m} t_k^{p/2} \mu(D_k) \\
	& \, \lesssim \,  \E \big ( \mathrm {W} _p^p( \mu_n^{R,T}, \mu ) \big )
 		+ \frac {1}{n^{p/d}}
\end {split} \eeq
where we used that $t_k = n^{-2/d} e^{k/d}$ and $\mu(D_k) \lesssim \, k^{(d/2)-1}e^{-k/2}$,
$k = 1, \ldots, m$.

As in the previous section, we thus concentrate on the study of
$\E ( \mathrm{ W}_p^p( \mu^{R, T}_n, \mu ) )$ that we control from
Proposition~\ref {prop.sobolev} together with the Riesz transform bound
\eqref {eq.sobolevriesz}. To this task, set
$$
g \, = \, g(y) \, = \,   f^{R, T}_n(y) - 1
	\, = \, \frac{1}{n} \sum_{i=1}^{n} \big [p_{T(X_i^R)} (X^R_i, y) - 1\big],
		\quad y \in \R^d.
$$
Therefore
$$
\E \big ( \mathrm {W} _p^p( \mu^{R, T}_n, \mu ) \big)
\, \leq \,  p^p \, \E \Big ( \big \| f^{R, T}_n - 1 \big \|^p_{H^{-1, p}(\mu)} \Big)
\, \lesssim \,  \E \bigg ( \int_{\R^d} \big | (- \mathrm {L})^{-1/2} g \big|^p d\mu \bigg).
$$
Center then the terms $p_{T(X_i^R)} (X^R_i, y)$ in the definition of
$g$ with respect to randomness in the $X_i^R$'s. To this task, write for every $y$,
$$
g(y) \, = \, \frac 1n\sum_{i=1}^{n} \big [p_{T(X_i^R)} (X^R_i, y)
		- \E (p_{T(X_i^R)} (X^R_i, y) ) \big]
+ \E \big (p_T(X^R, y) \big) - 1  \, = \,  \widetilde{g}(y) + \phi(y)
$$
with $\phi (y) = \E  (p_T(X^R, y)) - 1$ so that
$$
\E \big ( \mathrm {W} _p^p( \mu^{R, T}_n, \mu ) \big)
 \, \lesssim \, \E \bigg ( \int_{\R^d} \big | (- \mathrm {L})^{-1/2} \widetilde {g} \big|^p d\mu \bigg)
 		+ \int_{\R^d} \big | (- \mathrm {L})^{-1/2} \phi \big|^p d\mu .
$$

Rosenthal's inequality \cite {R70} for independent centered random variables $V_1, \ldots , V_n$
with a $p$-th moment, $p \geq 2$, expresses that
\beq \label {eq.rosenthal}
\E \bigg ( \bigg | \sum_{i=1}^n V_i \bigg |^p \bigg)
	\, \leq \, C_p \sum_{i=1}^n \E \big ( | V_i|^p \big)
			+ C_p \bigg ( \sum_{i=1}^n \E (  V_i^2) \bigg)^{p/2}
\eeq
where $C_p >0$ only depends on $p$.
For each fixed $y \in \R^d$, apply this inequality to the independent identically
distributed and centered random variables
$$
(- \mathrm {L}_y)^{-1/2} \big [p_{T(X_i^R)} (X^R_i, y)
		- \E (p_{T(X_i^R)} (X^R_i, y) ) \big], \quad i=1, \ldots, n,
$$
to get that
\beqs \begin {split}
\E \Big (\big | (- \mathrm {L}_y)^{-1/2} \, \widetilde {g}(y) \big|^p \Big)
  & \, \lesssim \, \frac {1}{n^{p-1}} \,
  	\E\Big ( \big | (- \mathrm {L}_y)^{-1/2} \big [p_T (X^R, y) - \E (p_T(X^R, y) ) \big] \big |^p \Big) \\
  & \quad \, \, +   \frac {1}{n^{p/2}} \,
  \E\Big ( \big [ (- \mathrm {L}_y)^{-1/2} \big [p_T (X^R, y) - \E(p_T(X^R, y)) \big] \big ]^2 \Big)^{p/2}.
\end {split} \eeqs
Since $\E (p_T(X^R, y) = \phi (y) +1$,
and we eventually aim to control $\int_{\R^d}  | (- \mathrm {L})^{-1/2} \phi |^p d\mu$, we may
replace back $\E (p_T(X^R, y) )$ by $1$ in the preceding. That is, we have at that point
\beq \begin {split} \label {eq.3termdecompositionp}
\E \big ( \mathrm {W} _p^p( \mu^{R, T}_n, \mu ) \big)
 & \, \lesssim \, \frac {1}{n^{p-1}} \int_{\R^d}
  	\E\Big ( \big | (- \mathrm {L}_y)^{-1/2} \big [p_T (X^R, y) - 1
		\big] \big |^p \Big) d\mu (y) \\
	& \quad \, + \frac {1}{n^{p/2}} \int_{\R^d}
  \E\Big ( \big [ (- \mathrm {L}_y)^{-1/2} \big [p_T (X^R, y) - 1 \big]
  		\big ]^2 \Big)^{p/2} d\mu (y) \\
 	& \quad \,	+ \int_{\R^d} \big | (- \mathrm {L})^{-1/2} \phi \big|^p d\mu .
\end {split} \eeq
In this expression, the random variables $(- \mathrm {L}_y)^{-1/2}  [p_T (X^R, y) - 1 ]$
will be studied with the help of the spectral representation \eqref {eq.spectral}, that is
$$
(- \mathrm {L}_y)^{-1/2} \big [p_T (X^R, y) - 1 \big]
	\, = \,  \frac {1}{\sqrt \pi}
	\int_0^\infty \! \frac {1}{\sqrt s} \, \big [p_{s+T} (X^R, y) - 1\big] ds
$$
and similarly for $(- \mathrm {L})^{-1/2} \phi $.

According to \eqref {eq.localizationregularizationp} and
\eqref {eq.3termdecompositionp}, the proof of the theorem will therefore be achieved once
it may be established that
\beq \label {eq.firsttermp}
\frac {1}{n^{p-1}}\int_{\R^d}
\E\bigg ( \bigg | \int_0^\infty  \! \frac {1}{\sqrt s} \,  \big [p_{s+T}(X^R, y) - 1 \big]
		ds\bigg |^p \bigg) d\mu(y) \, \lesssim \, \frac {1}{n^{p/d}} \, ,
\eeq
\beq \label {eq.secondtermp}
\frac {1}{n^{p/2}} \int_{\R^d}
\E\bigg ( \bigg [ \int_0^\infty \! \frac {1}{\sqrt s} \,  \big [p_{s+T}(X^R, y) - 1 \big]
		ds\bigg ]^2 \bigg)^{p/2} d\mu(y) \, \lesssim \, \frac {1}{n^{p/d}}
\eeq
and
\beq \label {eq.centeringp}
\int_{\R^d} \bigg | \int_0^\infty \frac {1}{\sqrt s} \, P_s \phi \, ds \bigg|^p d\mu
\, \lesssim \, \frac {1}{n^{p/d}} \, .
\eeq
The centering term \eqref {eq.centeringp} will be examined at the end of the proof.
We concentrate on the first two terms, starting with
the investigation of \eqref {eq.secondtermp} which is the most delicate one.

\bigskip

\noindent \textbf {Study of \eqref {eq.secondtermp}.}
Fix $y \in \R^d$ to begin with; by definition of the map $T$,
\beqs \begin {split}
\E\bigg ( \bigg [ \int_0^\infty \! \frac {1}{\sqrt s} \,  \big [p_{s+T} & (X^R, y) - 1 \big]
		ds\bigg ]^2 \bigg) \\
		& \, = \, \frac {1}{\mu(B_R)} \sum_{k=1}^m \int_{D_k}
	\bigg [ \int_0^\infty \! \frac {1}{\sqrt s} \,  \big [p_{s+t_k}(x, y) - 1 \big]  ds\bigg ]^2 d\mu(x).
\end {split} \eeqs
Given $s_k \in (0,1)$, $k = 1, \ldots, m$, to be specified, we decompose the integral in $s$
on $(0,s_k)$ and $(s_k, \infty)$ and study separately, by the triangle inequality, the resulting two pieces
in \eqref {eq.secondtermp}, showing that
\beq \label {eq.secondtermp<sk}
\frac {1}{n^{p/2}} \int_{\R^d}
\bigg ( \sum_{k=1}^m \int_{D_k}
	\bigg [ \int_0^{s_k} \! \frac {1}{\sqrt s} \,  \big [p_{s+t_k}(x, y) - 1 \big]  ds\bigg ]^2 d\mu(x)
	\bigg)^{p/2} d\mu(y) \, \lesssim \, \frac {1}{n^{p/d}}
\eeq
and
\beq \label {eq.secondtermp>sk}
\frac {1}{n^{p/2}} \int_{\R^d}
\bigg ( \sum_{k=1}^m \int_{D_k}
	\bigg [ \int_{s_k}^\infty \! \frac {1}{\sqrt s} \,  \big [p_{s+t_k}(x, y) - 1 \big]  ds\bigg ]^2 d\mu(x)
	\bigg)^{p/2} d\mu(y) \, \lesssim \, \frac {1}{n^{p/d}} \, .
\eeq

Concerning \eqref {eq.secondtermp<sk}, for each $k= 1, \ldots, m$ and $y \in \R^d$, by Fubini's theorem,
\beqs \begin {split}
\int_{D_k} \bigg [ \int_0^{s_k} \! & \frac {1}{\sqrt s} \,  \big [p_{s+t_k}(x, y) - 1 \big]  ds\bigg ]^2 d\mu(x)\\
	&\, = \, \int_{D_k}
	\int_0^{s_k} \! \int_0^{s_k} \frac {1}{\sqrt {ss'}}
		\,  \big [p_{s+t_k}(x, y) - 1 \big]   \big [p_{s'+t_k}(x, y) - 1 \big]  ds ds' d\mu(x) \\
	&\, \lesssim \, \int_0^{s_k} \! \int_0^{s_k} \frac {1}{\sqrt {ss'}}
		\int_{D_k} p_{s+t_k}(x, y) \, p_{s'+t_k}(x, y)  d\mu(x) \, ds ds'  + \mu (D_k). \\
\end {split} \eeqs
Summing over $k$, it is clear that the contribution $\mu(D_k)$ will be irrelevant for the final bound
and it is therefore ignored below. Now, a standard calculation on the explicit expression
of the Mehler kernel $p_t(x,y)$ yields
\beq \label {eq.extendedsemigroup}
\int_{D_k} p_{s+t_k}(x, y) \,  p_{s'+t_k}(x, y)  d\mu(x) \, = \,
    p_{s +s'+2t_k}(y, y) \,  \mu (\widetilde {D}_k)
\eeq
where $\widetilde {D}_k \, = \, - \frac \beta \alpha \, y + \alpha D_k$,
$$
\alpha^2  \, = \,  1 + \frac {a^2}{1 -a^2}  + \frac {b^2}{1 - b^2} \, , \qquad
\beta \, = \, \frac {a}{1 -a^2}  + \frac {b}{1 - b^2} \, ,
$$
with $a = e^{-s-t_k}$, $b = e^{-s'-t_k}$. For the further purposes, note that
\beq \label {eq.alphabeta}
\frac {\alpha^2}{\beta} - 1 \, = \, \frac {(1-a)(1-b)}{a+b} \, .
\eeq


We examine separately the expression $p_{s +s'+2t_k}(y, y)  \mu (\widetilde {D}_k) $
in \eqref {eq.extendedsemigroup} according
as $y \in \frac{\alpha^2}{\beta} D_k$ or not via the bound
\beq \label {eq.decomposition}
p_{s +s'+2t_k}(y, y) \,  \mu (\widetilde {D}_k)
	\, \leq \,	p_{s +s'+2t_k}(y, y) \,
	 \big [ \mathbbm {1}_{\frac{\alpha^2}{\beta} D_k} (y)
	+  \mathbbm {1}_{(\frac{\alpha^2}{\beta} D_k)^c} (y)  \mu (\widetilde {D}_k)   \big].
\eeq
As such, the study of \eqref {eq.secondtermp<sk} is divided into two parts, and the task is
to show that
\beq \label {eq.decompositionp1}
  \frac {1}{n^{p/2}} \int_{\R^d} \bigg [ \sum_{k=1}^m
		\int_0^{s_k} \! \int_0^{s_k} \!  \frac {1}{\sqrt {ss'}} \,
	 p_{s +s'+2t_k}( y, y)  \mathbbm {1}_{\frac{\alpha^2}{\beta} D_k} (y) ds ds'\bigg ]^{p/2} d\mu(y)
	 \, \lesssim \, \frac {1}{n^{p/d}}
\eeq
and
\beq \label {eq.decompositionp2}
 \frac {1}{n^{p/2}}  \int_{\R^d} \bigg [ \sum_{k=1}^m
		\int_0^{s_k} \! \int_0^{s_k} \! \frac {1}{\sqrt {ss'}} \,
	 p_{s +s'+2t_k}( y, y)  \mathbbm {1}_{(\frac{\alpha^2}{\beta} D_k)^c} (y)
	 \mu (\widetilde {D}_k)  ds ds' \bigg ]^{p/2} d\mu(y)
		\, \lesssim \, \frac {1}{n^{p/d}} \, .
\eeq

Start with \eqref {eq.decompositionp1} and fix now
$s_k = \frac {1}{\sqrt k}$, $k = 1, \ldots, m$.
By \eqref {eq.alphabeta}, it may be verified that for every $0 < s, s' \leq s_k$,
$$
\frac{\alpha^2}{\beta} D_k  \, \subset \,
      E_k \, = \; \big \{ x \in \R^d ; r_{k-1} \leq |x| < (1 + \textstyle {\frac 8k}) r_k \big \} .
$$
Hence, together with \eqref{eq.mehlerdiagonal} and a simple integration in $s,s'$, for every $y \in \R^d$,
$$
\int_0^{s_k} \! \int_0^{s_k} \! \frac {1}{\sqrt {ss'}} \,
	 p_{s +s'+2t_k}(y, y)  \mathbbm {1}_{\frac{\alpha^2}{\beta} D_k} (y) ds ds'
	\, \lesssim \, \mathbbm {1}_{ E_k} (y) \, \frac {1}{t_k^{(d/2)-1}} \, e^{|y|^2/2} .
$$
It follows that the left-hand side of \eqref {eq.decompositionp1} is bounded from above by
$$
	\frac {1}{n^{p/2}}  \int_{\R^d} \bigg [ \sum_{k=1}^m
		\mathbbm {1}_{ E_k} (y) \, \frac {1}{t_k^{(d/2)-1}}  \bigg]^{p/2} e^{p|y|^2/4} d\mu(y) .
$$
Observe that for every $k=1,\ldots , m$,
$ E_k \subset \bigcup_{\ell = k}^{k+80} D_\ell$
(with the obvious extension of $D_k$ when $k\geq m$).
As a consequence, the latter is bounded from above by
$$
\frac {1}{n^{p/2}}
	\sum_{k=1}^{m+80}  \frac {1}{t_k^{(\frac d2 -1)\frac p2}} \int_{D_k}  e^{p|y|^2/4} d\mu(y)
	\, \lesssim \, \frac {1}{n^{p/2}}
	\sum_{k=1}^{m+80}  \frac {1}{t_k^{(\frac d2 -1)\frac p2}}
			\, k^{\frac d2-1} e^{(\frac p2 - 1)\frac k2} .
$$
Now $t_k = n^{-2/d} e^{k/d}$, $k= 1, \ldots, m$, so that since $p < d$
the latter is of the order $\frac {1}{n^{p/d}}$, proving \eqref {eq.decompositionp1}.

We turn next to \eqref {eq.decompositionp2}. Fix $k = 1, \ldots, m$.
When $y \notin \frac{\alpha^2}{\beta} D_k$, then (draw a picture), it is clear that
$$
 \inf \big \{ |z| \in \R^d ; z \in \widetilde {D}_k \big \} 
 \, = \, \Big| \frac \beta\alpha \, |y| - \alpha \widetilde {r_k} \Big |
$$
where $ \widetilde {r_k} = r_k$ or $r_{k-1}$. Thus
$$
\mu (\widetilde {D}_k)
	 \, \lesssim \,
	 \bigg(1 + \Big| \frac \beta\alpha \, |y| - \alpha \widetilde {r_k} \Big|^{d-2} \bigg)
	 \exp \bigg ( - \frac 12 \Big ( \frac \beta\alpha \, |y| - \alpha \widetilde {r_k} \Big)^2 \bigg).
$$
To get rid of the prefactors in front of the exponential, 
let $\sigma \in (0, 1)$ (which will depend on $p$ and $d$ only) so that
$$
\mu (\widetilde {D}_k)
  \, \lesssim \,
   \exp \bigg ( - \frac {1-\sigma}{2} \Big ( \frac \beta\alpha \, |y| - \alpha \widetilde {r_k} \Big)^2 \bigg).
$$
Hence, together with \eqref {eq.mehlerdiagonal},
$$
p_{s+s'+2t_k}(y,y) \mathbbm {1}_{(\frac {\alpha^2}{\beta} D_k)^c} (y) \mu (\widetilde {D}_k)
  \, \lesssim \,
  \frac {1}{(1-a^2b^2)^{d/2}} \,
   \exp \bigg ( - \frac {1-\sigma}{2} \Big ( \frac \beta\alpha \, |y| - \alpha \widetilde {r_k} \Big)^2
   		+ \frac {ab}{1 +ab} \, |y|^2 \bigg).
$$

Now, since
$$
\frac {\beta^2}{\alpha^2} - \frac {2ab}{1 +ab} \, = \, \alpha^2 - 1,
$$
it is easily seen that
$$
 -\frac{1-\sigma}{2}\bigg (\frac{\beta}{\alpha}\, |y|^2 - \alpha \widetilde{r_k} \bigg)^2
 	 + \frac{ab}{1+ab} \, |y|^2
 	\, \leq  \, (1-\sigma)A + \frac{\sigma}{2} \, |y|^2
$$
where
$$
A \, = \, - \frac 12 (\alpha^2 - 1) |y|^2 + \beta \, \widetilde {r_k} |y|
		- \frac 12 \, \alpha^2 \widetilde {r_k}^2 .
$$
We would like to choose numerical non-negative constants $K$ and $L$ such that
$$
A \, \leq \, A' \, = \, \frac{K}{2}\, |y|^2 + \frac{L}{2} \, \widetilde{r_k}^2
$$
for all $y \in \R^d$ and $k=1, \ldots, m$. This may be achieved provided that $K + L \geq1$
and $L \leq 1$. Indeed, the quadratic form
$$
 B  \, = \, A' - A \, = \,  \frac{1}{2} (\alpha^2 - 1 + K) |y|^2
		- \beta \, \widetilde{r_k} |y| + \frac{1}{2} (\alpha^2+L) \widetilde{r_k}^2
$$
is positive semi-definite if $Q = \beta^2 - (\alpha^2 - 1 + K)(\alpha^2+L) \leq 0 $. But when $K + L \geq 1$ and $L \leq 1$,
\beqs \begin{split}
Q
 \, &= \,  - \frac{(1-a^2b^2)}{(1-a^2)(1-b^2)} \Big (K + L - \frac{2ab}{1+ab} \Big)  + L(1 - K)\\
 \, &\leq \, - \frac{(1-a^2b^2)}{(1-a^2)(1-b^2)} \Big (1 - \frac{2ab}{1+ab} \Big) + L^2\\
 \, &\leq \, - \frac{(1-ab)^2}{(1-a^2)(1-b^2)} + L^2 \leq -1 + L^2 \leq 0.
\end{split} \eeqs

As a consequence of this analysis, we have obtained that, provided $K + L  \geq 1$ and $L \leq 1$,
$$
p_{s+s'+2t_k} (y,y) \mathbbm {1}_{(\frac {\alpha^2}{\beta} D_k)^c} (y) \mu (\widetilde {D}_k)
  \, \lesssim \, \frac {1}{(1-a^2b^2)^{d/2}} \,  e^{A''}
  \, \lesssim \, \frac {1}{(s+s'+2t_k)^{d/2}} \,  e^{A''}
$$
uniformly over $s,s' \leq s_k$ and $y \in \R^d$, where
$$
A''  \, = \,  (1-\sigma)A' + \frac \sigma 2 |y|^2  \, = \,  \frac{(1 - \sigma)K + \sigma}{2} \, |y|^2
			+ \frac{(1 - \sigma)L}{2} \, \widetilde{r_k}^2 .
$$

It may now be integrated in $0 < s,s' \leq s_k$
for every $k$ to get that the left-hand side of \eqref {eq.decompositionp2} is bounded from above by
\beq \label {eq.Asecond}
 \frac {1}{n^{p/2}}
    \int_{\R^d}  \bigg[  \sum_{k=1}^m \frac {1}{t_k^{(d/2)-1}} \, e^{A''}\bigg ]^{p/2} d\mu(y)
\eeq
still under the condition $K + L \geq 1$ and $L \leq 1$.  The definition of $A''$ yields after integration that
\beqs \begin{split}
 \frac {1}{n^{p/2}}
    \int_{\R^d}  \bigg[  \sum_{k=1}^m \frac {1}{t_k^{(d/2)-1}} & \, e^{A''}\bigg ]^{p/2} d\mu(y) \\
 	\, &= \,   \frac {1}{n^{p/2}} \int_{\R^d} \bigg[  \sum_{k=1}^m \frac {1}{t_k^{(d/2)-1}} \, 				e^{\frac{(1-\sigma)L}{2}\widetilde{r_k}^2} \bigg]^{p/2}
					e^{\frac{p[(1-\sigma)K+\sigma]}{4} |y|^2} d\mu(y)\\
 	\, &\lesssim \, \frac {1}{n^{p/2}} \frac {1}{(1 - \frac {p}{2}[(1-\sigma)K + \sigma ])^{d/2} } 					\, \bigg[  \sum_{k=1}^m \frac {1}{t_k^{(d/2)-1}}
 					\, e^{\frac{(1-\sigma)L}{2}\widetilde{r_k}^2} \bigg]^{p/2}
\end{split} \eeqs
provided that $\frac {p}{2} [(1-\sigma)K + \sigma ]< 1$.
Since $t_k = n^{-2/d}e^{k/d}$, $k=1, \ldots, m$, and  $\widetilde {r_k} = \sqrt {k} $ or $\sqrt {k-1}$,
it is necessary in addition that $ (1 - \sigma)L < 1 - \frac 2d$ in order that
the preceding bound yields the correct rate $\frac {1}{n^{p/d}}$. Provided that  $\sigma \in (0, \frac{2}{p})$, the
two preceding conditions are indeed compatible with $K + L \geq 1$ and $L \leq 1$. As an interesting example, we can always set $K = 0$, $L = 1$ and $\sigma \in (\frac{2}{d}, \frac{2}{p})$,
simplifying therefore the exposition.
The preceding construction will however be needed later in the case $p=d$ (Section~\ref {sec.6}).
As a conclusion, \eqref {eq.decompositionp2} is established.

As announced, the two controls \eqref {eq.decompositionp1} and \eqref {eq.decompositionp2}
yield together the expected bound \eqref {eq.secondtermp<sk}.

\medskip

We turn to the analysis of the second part \eqref {eq.secondtermp>sk} concerned with
the values of $s \geq s_k$. The first step is a Minkowski integral inequality to
exchange the order of integration. To this purpose, it is convenient to rewrite
the left-hand side of \eqref {eq.secondtermp>sk} as
\beq \label {eq.secondtermp>skbis}
\frac {1}{n^{p/2}} \int_{\R^d}
   \bigg( \E \bigg( \bigg | \int_S^{\infty} \frac {1}{\sqrt s} \,  \big [p_{s+T}(X^R, y) - 1 \big]ds \bigg |^2
    \bigg)    \bigg)^{p/2} d\mu(y)
\eeq
where $S  : B_R \to (0,1)$ is defined by $S(x) = s_k$ if $x \in D_k$, $k = 1, \ldots, m$.
As for $T=T(X^R)$, we write here $S = S(X^R)$ to ease the notation. As announced, by Minkowski's inequality
since $p \geq 2$, the latter is less than or equal to
$$
\frac {1}{n^{p/2}}  \,  \E \bigg (  \bigg( \int_{\R^d} \bigg | \int_S^{\infty} \frac {1}{\sqrt s} \,
		 \big [p_{s+T}(X^R, y) - 1 \big]ds \bigg |^p d\mu(y) \bigg)^{2/p} \bigg)^{p/2}.
$$
Now, and conditionally on the randomness of $X^R$, also by the triangle inequality,
$$
\int_{\R^d} \bigg | \int_S^{\infty} \frac {1}{\sqrt s} \,
		 \big [p_{s+T}(X^R, y) - 1 \big]ds \bigg |^p d\mu(y)
		 \, \leq \, \bigg ( \int_S^\infty \frac {1}{\sqrt s} \,
		 {\big \| p_{s+T}(X^R, \cdot) - 1 \big \|}_p \, ds \bigg)^p
$$
where the norm ${\| \cdot \|}_p$ is in $\mathrm {L}^p(d\mu(y))$.
By \eqref {eq.hypercontractivityp}, for every $s \geq S$,
$$
{\big \| p_{s+T}(X^R, \cdot) - 1 \big \|}_p
	\, \lesssim \, e^{-s/2} {\big \| p_S(X^R, \cdot) - 1 \big \|}_p
$$
so that
\beq \label {eq.hypercontractivityrandomp}
\bigg ( \int_S^\infty \frac {1}{\sqrt s} \,
		 {\big \| p_{s+T}(X^R, \cdot) - 1 \big \|}_p \, ds \bigg)^p
		 \, \lesssim \, {\big \| p_S(X^R, \cdot) - 1 \big \|}_p^p
		 \, \lesssim \, \frac {1}{S^{(p-1)d/2}} \, e^{(p-1)|X^R|^2/2}
\eeq
where \eqref {eq.mehlerintegral} is used in the last step.
Therefore, with the choice of $s_k = \frac {1}{\sqrt k}$,
\eqref {eq.secondtermp>skbis} is bounded from above by
\beqs \begin {split}
\frac {1}{n^{p/2}} \, \E \bigg( \frac{1}{S^{(1-\frac 1p)d}} \, e^{(1-\frac{1}{p})|X^R|^2}  \bigg)^{p/2}
	&  \, \lesssim \, \frac {1}{n^{p/2}}
	 \bigg ( \sum_{k=1}^m k^{(1-\frac 1p)\frac d2}
	 	 \int_{D_k} e^{(1-\frac{1}{p})|x|^2} d\mu(x) \bigg)^{p/2} \\
 	& \, \lesssim \, \frac {1}{n^{p/2}}
 \bigg( \sum_{k=1}^{m} k^{(1 - \frac {1}{2p})d - 1} \, e^{(\frac{1}{2} - \frac{1}{p})k } \bigg)^{p/2}.
\end {split} \eeqs
Recalling that $m =R^2 = 2 c \log n$ with $ c <1$ yields an expression of at
most the order $\frac {1}{n^{p/d}}$ so that \eqref {eq.secondtermp>sk} is established.

As a consequence of \eqref {eq.secondtermp<sk} and \eqref {eq.secondtermp>sk},
the bound \eqref {eq.secondtermp} is established.

\bigskip

\noindent \textbf {Study of \eqref {eq.firsttermp}.}
We address here \eqref {eq.firsttermp} following the steps developed
for \eqref {eq.secondtermp} but in a simplified way.
Decomposing the integral in $s$, simply here on $(0,1)$ and $(1, \infty)$
as in the case $p=2$, \eqref {eq.firsttermp} will hold as soon as
\beq \label {eq.firsttermp<1}
\frac {1}{n^{p-1}} \, \E \bigg ( \int_{\R^d} \bigg | \int_0^1 \frac {1}{\sqrt s} \,
	\big [ p_{s+T} (X^R, y) - 1 \big] \bigg |^p d\mu (y) \bigg)
	 \, \lesssim \, \frac {1}{n^{p/d}}
\eeq
and
\beq \label {eq.firsttermp>1}
\frac {1}{n^{p-1}} \, \E \bigg ( \int_{\R^d} \bigg | \int_1^\infty \frac {1}{\sqrt s} \,
	\big [ p_{s+T} (X^R, y) - 1 \big] \bigg |^p d\mu (y) \bigg)
\, \lesssim \, \frac {1}{n^{p/d}} \, .
\eeq
The main simplification here with respect to \eqref {eq.secondtermp} is that
Fubini's theorem applies between the expectation $\E$ and integration in $d\mu(y)$.

Starting with \eqref {eq.firsttermp<1}, recall that by definition of the map $T$,
for every $y \in \R^d$,
\beqs \begin {split}
\E\bigg ( \bigg |\int_0^1 \frac {1}{\sqrt s} \,  \big [p_{s+T} & (X^R, y) - 1 \big]
		ds\bigg |^p \bigg) \\
		& \, = \, \frac {1}{\mu(B_R)} \sum_{k=1}^m \int_{D_k}
	\bigg | \int_0^1 \frac {1}{\sqrt s} \,  \big [p_{s+t_k}(x, y) - 1 \big]  ds\bigg |^p d\mu(x).
\end {split} \eeqs
Then, by the triangle inequality, for every $k= 1, \ldots, m$,
$$
\int_{D_k} \bigg | \int_0^1 \frac {1}{\sqrt s} \,  \big [p_{s+t_k}(x, y) - 1 \big]  ds\bigg |^p d\mu(x)
	\, \lesssim \, \int_{D_k}
	\bigg [ \int_0^1 \frac {1}{\sqrt s} \, p_{s+t_k}(x, y)  ds\bigg ]^p d\mu(x)
	+ \mu(D_k)
$$
and it is clear again that we may ignore the contribution $\mu(D_k)$ in what follows.
By \eqref {eq.mehleruniform}, for all $x, y \in \R^d$,
$$
p_{s+t_k}(x, y) \, \leq \,  \frac {1}{(1 - e^{-2(s+t_k)})^{d/2}} \, e^{|x|^2/2}
			\, \lesssim \,  \frac {1}{(s+t_k)^{d/2}} \, e^{|x|^2/2}
$$
in the range $s \leq 1$. Hence, after integration in $s$, for every $x,y \in \R^d$,
$$
	 \int_0^1 \frac {1}{\sqrt s} \, p_{s+t_k}(x, y)  ds
 		 \, \lesssim \,  \frac {1}{t_k^{(d-1)/2}} \, e^{|x|^2/2}.
$$
Therefore, for every $k$,
\beqs \begin {split}
\int_{D_k}  \int_{\R^d}
	\bigg [ \int_0^1 \frac {1}{\sqrt s} \, &  p_{s+t_k}(x, y)   ds\bigg ]^p d\mu(y)
	d\mu(x) \\
	 & \, \lesssim \,  \frac {1}{t_k^{(d-1)(p-2)/2}}
	\int_{D_k}  \int_{\R^d}
	\bigg [ \int_0^1 \frac {1}{\sqrt s} \,  p_{s+t_k}(x, y)   ds\bigg ]^2 d\mu(y)
	 \, e^{(p-2)|x|^2/2}  d\mu(x).
\end {split} \eeqs
Now
$$
\int_{\R^d}
	\bigg [ \int_0^1 \frac {1}{\sqrt s} \,  p_{s+t_k}(x, y)   ds\bigg ]^2 d\mu(y)
	\, = \, \int_0^1 \! \int_0^1 \frac {1}{\sqrt {ss'}} \, p_{s+s'+2t_k} (x,x) ds ds'
$$
is, by \eqref {eq.mehlerdiagonal}, of the order of at most $\frac {1}{t_k^{(d/2)-1}} \,e^{|x|^2/2}$.
It follows that
\beqs \begin {split}
\frac {1}{n^{p-1}}
 \sum_{k=1}^m \int_{D_k}  \int_{\R^d}
	\bigg [ \int_0^1 \frac {1}{\sqrt s} \,  & p_{s+t_k}(x, y)   ds\bigg ]^p d\mu(y) d\mu(x) \\
   	&\, \lesssim \, \frac {1}{n^{p-1}} \sum_{k=1}^m \frac {1}{t_k^{(pd -p -d)/2}}
		\int_{D_k} e^{(p-1) |x|^2/2} d\mu(x) \\
 	& \, \lesssim \, \frac {1}{n^{p-1}} \sum_{k=1}^m
 			\frac {k^{(d/2)-1}}{t_k^{(pd -p -d)/2}}\,  e^{(p-2) k/2} \\
  	& \, \lesssim \, \frac {1}{n^{p/d}} \\
\end {split} \eeqs
by the choice of $t_k = n^{-2/d} e^{k/d}$, $k=1, \ldots, m$, together with the fact
that $p < d$. Hence \eqref {eq.firsttermp<1} holds true.

Concerning \eqref {eq.firsttermp>1}, the arguments developed for \eqref {eq.secondtermp>skbis}
may essentially be repeated.
In particular, making use of \eqref {eq.hypercontractivityrandomp}, the left-hand
side of \eqref {eq.firsttermp>1} may be seen to be bounded from above by
$$
\frac {1}{n^{p-1}} \, \E \big(  e^{(p-1)|X^R|^2/2}  \big)
	  \, \lesssim \, \frac {1}{n^{p-1}}
	 \sum_{k=1}^m  \int_{D_k} e^{(p-1)|x|^2/2} d\mu(x)
 	 \, \lesssim \, \frac {1}{n^{p-1}}  \sum_{k=1}^{m}  k^{(d/2)-1} e^{(p-2)k/2 }    .
$$
Again since $m =R^2 = 2 c \log n$ with $ c <1$, this contribution is at
most the order $\frac {1}{n^{p/d}}$ proving \eqref {eq.firsttermp>1}.

As announced, as a consequence of \eqref {eq.firsttermp<1} and \eqref {eq.firsttermp>1},
the bound \eqref {eq.firsttermp} is established.

\bigskip

\noindent \textbf {Study of \eqref {eq.centeringp}.}
In the final part of the proof, we thus take care of the centering term
$$
\int_{\R^d} \bigg | \int_0^\infty \frac {1}{\sqrt s} \, P_s \phi \, ds \bigg|^p d\mu
$$
of \eqref {eq.centeringp}. Recall that by definition
$$
\phi \, = \, \E \big ( p_T(X^R, \cdot) \big) - 1
 	\, = \,  \frac{1}{\mu(B_R)} \sum_{k = 1}^{m} \big [P_{t_k}(\mathbbm {1}_{D_k}) - \mu(D_k) \big],
$$
and write then
\beqs \begin {split}
P_s \phi
	& \,=\,  \frac{1}{\mu(B_R)} \sum_{k = 1}^{m} \big [P_{s+t_k}(\mathbbm {1}_{D_k} )
 - P_s (\mathbbm {1}_{D_k}) \big ]
  + \frac{1}{\mu(B_R)} \big [ P_s (\mathbbm {1}_{B_R}) - \mu(B_R) \big ] \\
  	& \, = \, \frac{1}{\mu(B_R)} \, \big ( \phi_{s,1} + \phi_{s,2} \big) .
\end {split} \eeqs
By means of H\"older's inequality (in $ds$), for any $0 < \kappa < \frac p2 - 1$,
\beq  \begin {split} \label {eq.holderkappa}
\bigg |  \int_0^\infty \frac {1}{\sqrt s} \, P_s \phi \, ds \bigg |^p
	&\, \lesssim \,  \int_0^\infty  s^{\frac p2-1 -\kappa} \, e^{\kappa s} \, | P_s \phi |^p  ds \\
	&\, \lesssim \,  \int_0^\infty  s^{\frac p2-1 -\kappa} \, e^{\kappa s} \, | \phi_{s,1} |^p  ds
			+  \int_0^\infty  s^{\frac p2-1 -\kappa} \, e^{\kappa s} \, | \phi_{s,2} |^p  ds .
\end {split} \eeq

We examine successively the contributions of $\phi_{s,1}$ and $ \phi_{s,2}$ in the preceding.
By the triangle inequality,
\beqs \begin {split}
\int_{\R^d} | \phi_{s,1}|^p d\mu
	& \, = \, \int_{\R^d} \bigg | \sum_{k = 1}^{m} \big [P_{s+t_k}(\mathbbm {1}_{D_k} )
 	- P_s (\mathbbm {1}_{D_k}) \big ] \bigg |^p d\mu \\
 	& \, \leq \, \bigg ( \sum_{k=1}^m { \big \|  P_{s+t_k}(\mathbbm {1}_{D_k} )
 		- P_s (\mathbbm {1}_{D_k}) \big \|}_p   \bigg )^p .\\
\end {split} \eeqs
Using \eqref {eq.pseudopoincare1}, for any $s>0$ and since $t_k < 1$,
$$
{ \big \|  P_{s+t_k}(\mathbbm {1}_{D_k} )
 		- P_s (\mathbbm {1}_{D_k}) \big \|}_p
		\, \lesssim \, \frac {e^{-s/p}}{(e^{2s} - 1)^{(p-1)/2p}} \, \sqrt{t_k}
		\, \mu(\partial D_k)^{1/p} .
$$
For $\kappa >0$ small enough, it follows that
\beq \label {eq.phip1}
\int_0^\infty s^{\frac p2-1 -\kappa} \, e^{\kappa s}  \int_{\R^d} | \phi_{s,1} |^p d\mu \,  ds
 \, \lesssim \, \bigg ( \sum_{k=1}^m \sqrt { t_k} \, \mu(\partial D_k)^{1/p} \bigg)^p
 \, \lesssim \, \frac {1}{n^{p/d}}
\eeq
since $t_k = n^{-2/d}e^{k/d}$ and $\mu( \partial D_k) \lesssim  k^{(d-1)/2} e^{-k/2}$,
$k=1, \ldots, m$.

On the other hand, it is easily seen as in the case $p=2$ that, again for $\kappa >0$ small enough,
\beq \begin {split} \label {eq.phip2}
 \int_0^\infty  s^{\frac p2-1 -\kappa} \, e^{\kappa s}
 		 \int_{\R^d} | \phi_{2,s} |^p d\mu \,  ds
	& \, \lesssim \,  \int_0^\infty  s^{\frac p2-1 -\kappa} \, e^{\kappa s}
 		 \int_{\R^d} | \phi_{2,s} |^2 d\mu \,  ds \\
	& \, \lesssim \,  \int_0^\infty  s^{\frac p2-1 -\kappa} \, e^{\kappa s- 2s}
 		 \, \mu (B_R^c)  ds  \\
	& \, \lesssim \, \frac {1}{n^{p/d}} \\
\end {split} \eeq
by the choice of $R = \sqrt {2c\log n}$.

Together with \eqref {eq.holderkappa}, it follows from \eqref {eq.phip1} and \eqref {eq.phip2} that
the bound \eqref {eq.centeringp} is established. Altogether, the proof of Theorem~\ref {thm.main}
is complete.

\section{The case $p = d$} \label {sec.6}
\setcounter{equation}{0}

This section addresses the proof of Theorem~\ref {thm.pd} for, thus, $p=d$ ($\geq 2$).
The proof for $p=d=2$ was actually provided in \cite {L17}, as a simpler version
of what is developed here. (In the next section,
we present the pde-transportation argument for the lower bound in this case.)
The proof here for $p=d$ carefully adjusts several parameters in the various
steps of the one developed in Section~\ref{sec.5}.

The first step is truncation on a ball $B_R$ this time with $R = \sqrt{ 2\log n}$ for which
it holds similarly that
\beq  \label{eq.localizationd}
\E \big ( \mathrm {W}_d^d( \mu_n, \mu^R_n ) \big)
				\, \lesssim \, \frac {( \log n)^{\frac{d}{2}}}{n} \, .
\eeq
The ball $B_R$ is decomposed again as the union of $m$ annuli
$D_k = \{x \in \R^d ; r_{k-1} \leq |x| < r_k \}$,
$k = 1, \ldots, m$, where $0 = r_0 < r_1 < \cdots <r_m = R$. Consider as well
the map $T : B_R \to (0, 1)$ defined by $T(x) = t_k$ if $x \in D_k$. We will use
the same choices $r_k = \sqrt {k}$, but modify the values of $t_k$ as
$$
t_k \, = \,  \frac{e^{k/d}}{ n^{2/d} \sqrt{k}} \, , \quad k = 1, \ldots, m.
$$

Setting
$$
f(y) \, = \,  f^{R, T}_n(y) \, = \,  \frac{1}{n} \sum_{i=1}^{n} p_{T(X^R_i)} (X^R_i, y), \quad y \in \R^d,
$$
and $d\mu^{R, T}_n = f^{R, T}_n d\mu$, the preceding choice of the $t_k$'s now yields
in \eqref {eq.regularizationp} that
\beq \label {eq.regularizationd}
\E \big ( \mathrm {W} _d^d( \mu^R_n, \mu^{R, T}_n ) \big )
 	  \, \lesssim \, \frac {(\log n)^{\max ( \frac{d}{2}-1, \frac{d}{4} )}}{n} \, .
\eeq

Next, to estimate $\E ( ( \mathrm {W} _d^d( \mu^{R, T}_n, \mu ) )$
as in Section~\ref{sec.5}, we need to control the terms
\beq \label {eq.firsttermd}
\frac {1}{n^{d-1}}\int_{\R^d}
\E\bigg ( \bigg | \int_0^\infty  \! \frac {1}{\sqrt s} \,  \big [p_{s+T}(X^R, y) - 1 \big]
		ds\bigg |^d \bigg) d\mu(y)  \, ,
\eeq
\beq \label {eq.secondtermd}
\frac {1}{n^{d/2}} \int_{\R^d}
\E\bigg ( \bigg [ \int_0^\infty \! \frac {1}{\sqrt s} \,  \big [p_{s+T}(X^R, y) - 1 \big]
		ds\bigg ]^2 \bigg)^{d/2} d\mu(y) \,
\eeq
and
\beq \label {eq.centeringd}
\int_{\R^d} \bigg | \int_0^\infty \frac {1}{\sqrt s} \, P_s \phi \, ds \bigg|^d d\mu
\,
\eeq
where $\phi(y) = \E \left (p_T(X^R, y) \right) - 1$, $ y \in \R^d$.

\bigskip

\noindent \textbf {Study of \eqref {eq.secondtermd}.} Given
$s_k = \frac {1}{\sqrt k}$, $k = 1, \ldots, m$, it is sufficient to estimate separately
\beq \label {eq.secondtermd<sk}
\frac {1}{n^{d/2}} \int_{\R^d}
\bigg ( \sum_{k=1}^m \int_{D_k}
	\bigg [ \int_0^{s_k} \! \frac {1}{\sqrt s} \,  \big [p_{s+t_k}(x, y) - 1 \big]  ds\bigg ]^2 d\mu(x)
	\bigg)^{d/2} d\mu(y) \,
\eeq
and
\beq \label {eq.secondtermd>sk}
\frac {1}{n^{d/2}} \int_{\R^d}
\bigg ( \sum_{k=1}^m \int_{D_k}
	\bigg [ \int_{s_k}^\infty \! \frac {1}{\sqrt s} \,  \big [p_{s+t_k}(x, y) - 1 \big]  ds\bigg ]^2 d\mu(x)
	\bigg)^{d/2} d\mu(y) \, .
\eeq

Concerning \eqref{eq.secondtermd<sk}, with the notation of the previous section, we need to investigate
\beq \label {eq.decompositiond1}
  \frac {1}{n^{d/2}} \int_{\R^d} \bigg [ \sum_{k=1}^m
		\int_0^{s_k} \! \int_0^{s_k} \!  \frac {1}{\sqrt {ss'}} \,
	 p_{s +s'+2t_k}( y, y)  \mathbbm {1}_{\frac{\alpha^2}{\beta} D_k} (y) ds ds'\bigg ]^{d/2} d\mu(y)
\eeq
and
\beq \label {eq.decompositiond2}
 \frac {1}{n^{d/2}}  \int_{\R^d} \bigg [ \sum_{k=1}^m
		\int_0^{s_k} \! \int_0^{s_k} \!  \frac {1}{\sqrt {ss'}} \,
	 p_{s +s'+2t_k}( y, y)  \mathbbm {1}_{(\frac{\alpha^2}{\beta} D_k)^c} (y)
	 	\mu (\widetilde {D}_k) ds ds'\bigg ]^{d/2} d\mu(y) \, .
\eeq

Arguing as for \eqref {eq.decompositionp1}, \eqref {eq.decompositiond1} is upper bounded by
$$
\frac{1}{n^{d/2}} \sum_{k=1}^{m+80}  \frac {1}{t_k^{(\frac d2 -1)\frac d2}}
			\, k^{\frac d2-1} e^{(\frac d2 - 1)\frac k2}
    \, \lesssim \, \frac{ (\log n)^{\frac{d^2 + 2d}{8}} }{n} \, .
$$
Turning to \eqref {eq.decompositiond2}, fix $k = 1, \ldots, m$. The proof proceeds as
for \eqref {eq.decompositionp2} in Section~\ref {sec.5} but now with a choice of
$K = K_k $ and $L = L_k $
this time depending on $k$, to reach that \eqref {eq.decompositiond2} is upper bounded by
$$
 \frac{1}{n}   \bigg [ \sum_{k=1}^m k^{\frac{d -2}{4}}
 	\frac {e^{[(1-\sigma)L_k - 1 +\frac{2}{d}]\frac{k}{2}}}{1 - \frac {d}{2}[(1-\sigma)K_k + \sigma ] }      	\bigg]^{d/2}
$$
(where we recall that $\sigma \in (0, \frac 2d)$).
Fix then $\varepsilon > 0$ small enough, and, for $k = 1, \ldots, m$, let
$$
L_k  \, = \, \frac{1 - \frac{2}{d} + \frac{\varepsilon}{k}}{1 - \sigma} \, , \quad
	K_k \, = \, 1 - L_k = \frac{ \frac{2}{d}- \sigma - \frac{\varepsilon}{k}}{1 - \sigma}
$$
so that $K_k + L_k = 1$ and $L_k < 1$. With these choices,
\beq
\frac{1}{n}   \bigg [ \sum_{k=1}^m k^{\frac{d -2}{4}}
	\frac {e^{[(1-\sigma)L_k - 1 +\frac{2}{d}]\frac{k}{2}}}{1 - \frac {d}{2}[(1-\sigma)K_k + \sigma ] }  		\bigg ]^{d/2}
	\, \lesssim \, \frac{1}{n} \bigg ( \sum_{k=1}^m k^{\frac{d + 2}{4}}  \bigg)^{d/2}
	\, \lesssim \, \frac{(\log n)^{\frac{d^2 + 6d }{8}}}{n} \, .
\eeq

As a consequence of the previous analysis, \eqref {eq.secondtermd<sk} is thus
controlled by $\frac{(\log n)^{\frac{d^2 + 6d }{8}}}{n}$. Concerning
\eqref {eq.secondtermd>sk}, it is handled as in Section~\ref {sec.5}
for the study of \eqref {eq.secondtermp>sk} and upper bounded by
$$
\frac {1}{n^{d/2}} \bigg( \sum_{k=1}^{m} k^{d-\frac 32}
		\, e^{(\frac{1}{2} - \frac{1}{d})k } \bigg)^{d/2}
		\, \lesssim \, \frac{(\log n)^{\frac {d^2}{2} - \frac {3d}{4}}}{n} \, .
$$
These two bounds lead to the dichotomy $d=2,3$ and $d \geq 4$ in the statement
of Theorem~\ref {thm.pd}.

\medskip

\noindent \textbf {Study of \eqref {eq.firsttermd}.} This term is simpler than \eqref {eq.secondtermd}. Following the analysis of \eqref {eq.firsttermp} in Section~\ref {sec.5}, the leading term is
$$
\frac {1}{n^{d-1}} \sum_{k=1}^m \frac {k^{(d/2)-1}}{t_k^{(d^2 - 2d)/2}}\,  e^{(d-2) k/2}
		 \, \lesssim \, \frac{(\log n)^{\frac{d^2}{4}}}{n}
$$
by the choice of $t_k = \frac{e^{k/d}}{ n^{2/d} \sqrt{k}} $, $k=1, \ldots, m$.

\medskip

\noindent \textbf {Study of \eqref {eq.centeringd}.} Using the method in the previous section,
it is bounded from above by
$$
\bigg ( \sum_{k=1}^m \sqrt { t_k} \, \mu(\partial D_k)^{1/d} \bigg)^d
			\, \lesssim \, \frac{(\log n)^{\frac{5d - 4}{4}}}{n} \, .
$$
since $t_k = \frac{e^{k/d}}{ n^{2/d} \sqrt{k}}$ and
$\mu( \partial D_k) \lesssim  k^{(d-1)/2} e^{-k/2}$, $k=1, \ldots, m$.

\medskip

To conclude, it may be checked that all the logarithmic exponents are less than or equal to
$$
\kappa = \max \Big ( \frac{d^2 + 6d }{8}  , \frac {d^2}{2} - \frac {3d}{4} \Big),
$$
thereby completing the proof of Theorem~\ref {thm.pd}.

\section {Lower bound $p=d=2$} \label {sec.7}

The purpose of this paragraph is to provide an
alternate proof based on the pde-transportation method of \cite {AST19} of the lower bound
\beq \label {eq.lowerboundd=2}
    \E  \big( \mathrm {W}_2^2 (\mu_n , \mu )\big )  \, \gtrsim \,  \frac {(\log n)^2}{n}
\eeq
of \eqref {eq.gaussian2} in dimension $2$, which has been established in \cite {T18}
by different means. This proof already appeared in \cite {L18}, and is included here for
completeness and convenience.

The first step is a two-sided bound on the Kantorovich metric $ \mathrm {W}_2$
in terms of Sobolev norms. It is developed in \cite {L18}
in weighted Riemannian manifolds
under the curvature condition $CD(K,\infty)$ for some $K \in \R$, but for simplicity
is restricted here to the Gaussian model (for which $K=1$).
Let thus $\mu$ be the standard Gaussian measure on the Borel sets of $\R^d$,
and $\mathrm {L}$ be the Ornstein-Uhlenbeck operator as presented in Section~\ref {sec.2}.
Proposition~\ref {prop.lowerbound} is taken from~\cite {AST19}.
Proposition~\ref {prop.upperbound} (a slight extension of Proposition~\ref {prop.sobolev})
is not used below but included for comparison.

\begin {proposition}  \label {prop.upperbound}
Let $d\nu = f d\mu$ and $f = 1 + g$, and let $0 < c \leq 1$. If $g \geq - c$, then
\beq \label {eq.upperbound}
\mathrm {W}_2^2 (\nu, \mu) \, \leq \,
  \frac {4}{c^2} \, \big[1 - \sqrt {1-c} \,\big ]^2 \int_{\R^d} g (- \mathrm{L} )^{-1} g \, d\mu
\eeq
(where $g$ is assumed to belong to the suitable domain so that the left-hand side makes sense).
\end {proposition}

\begin {proposition}  \label {prop.lowerbound}
Let $d\nu = f d\mu$ and $f = 1 + g$.
Then, whenever $g$ and $h : \R^d \to \R$ belong to the suitable domain and $h$ is such that
$\int_{\R^d} h d\mu = 0$ and $h \leq c$ uniformly for some $c >0$,
\beq \label {eq.lowerbound1}
 \mathrm {W}_2^2 (\nu, \mu)
   \, \geq \, 2 \int_{\R^d} g (- \mathrm{L} )^{-1} h \, d\mu
     -   \frac {e^{c} -1}{c}   \int_{\R^d}  h (- \mathrm{L} )^{-1} h \, d\mu .
\eeq
In particular, if $g \leq c$,
\beq \label {eq.lowerbound3}
 \mathrm {W}_2^2 (\nu, \mu)
    \, \geq \, \bigg ( 2   -   \frac {e^{c} -1}{c}  \bigg)  \int_{\R^d} g (- \mathrm{L} )^{-1} g \, d\mu .
\eeq
\end {proposition}

Recall that by integration by parts
$$
\int_{\R^d} g (- \mathrm{L} )^{-1} g \, d\mu \, = \,
		\int_{\R^d}  \big| \nabla ((-\mathrm {L})^{-1}g) \big|^2  d\mu
$$
which is the Sobolev norm of Proposition~\ref {prop.sobolev}.
Note also that as $c \to 0$,
$$
\frac {4}{c^2} \, \big[1 - \sqrt {1-c} \,\big ]^2 \, \sim \, 1 + \frac c2
   		\qquad \mbox {and} \qquad  2 -  \frac {e^c -1}{c} \, \sim \,  1 - \frac c2
$$
so that the bounds \eqref {eq.upperbound} and \eqref {eq.lowerbound3} are sharp in this regime.

\begin {proof} [Proof of Proposition~\ref {prop.upperbound}]
It is shown in \cite {L17} that for every (smooth) increasing
$\theta : [0,1] \to [0,1]$ with $\theta (0) = 0$, $\theta (1) = 1$,
$$
 \mathrm {W}_2^2 (\nu, \mu)
    \, \leq \,  \int_{\R^d}  \big| \nabla ((-\mathrm {L})^{-1}g) \big|^2
       \int_0^1 \frac {\theta'(s)^2}{1 + \theta(s)g} \, ds \, d\mu .
$$
Using that $g \geq -c$,
$$
 \mathrm {W}_2^2 (\nu, \mu)
    \, \leq \,   \int_0^1 \frac {\theta'(s)^2}{1 - \theta(s)c} \, ds
    \int_{\R^d}  \big| \nabla ((-\mathrm {L})^{-1}g) \big|^2   d\mu .
$$
The claim \eqref {eq.upperbound} follows from the (optimal) choice
$$
\theta (s) \, = \, \frac {1 - \sqrt {1-c}}{c} \, \Big ( 2s - \big[1 - \sqrt {1-c} \,\big ] s^2 \Big) ,
		\quad s \in [0,1].
$$
When $c=1$, the conclusion amounts to Proposition~\ref {prop.sobolev}.
\end {proof}

\begin {proof} [Proof of Proposition~\ref {prop.lowerbound}]
As announced, we follow \cite {AST19}.
By the Kantorovich dual description of the Kantorovich
metric $\mathrm {W}_2$ (cf.~\cite {V09}), for any bounded continuous $ \varphi : \R^d \to \R$,
\beqs \begin {split}
\frac 12 \, \mathrm {W}_2^2 (\nu, \mu)
		& \, \geq \, \int_{\R^d} \varphi \, f d\mu - \int_{\R^d} \widehat {Q}_1 \varphi \, d\mu \\
		& \, = \, \int_{\R^d} \varphi \, g \, d\mu
		  - \bigg( \int_{\R^d} \widehat {Q}_1 \varphi \, d\mu - \int_{\R^d} \varphi \, d\mu \bigg)
\end {split} \eeqs
where $\widehat {Q}_1$ is the supremum convolution
$$
\widehat {Q}_1 \varphi (x) \, = \, \sup_{y \in \R^d} \Big [ \varphi (y) - \frac 12 \, |x-y|^2 \Big].
$$
Choose then $ \varphi = (- \mathrm{L})^{-1} h$. Now
$$
\int_{\R^d} \widehat {Q}_1 \varphi \, d\mu - \int_{\R^d} \varphi \, d\mu
		\, = \, \frac 12 \int_0^1 \int_{\R^d} | \nabla \widehat {Q}_s \varphi |^2 d\mu \, ds.
$$
It is shown in \cite {AST19} that since $ -\mathrm {L} \varphi =   h \leq c$ uniformly,
under a $CD(0, \infty)$ curvature condition,
$$
\int_{\R^d} | \nabla \widehat {Q}_s \varphi |^2 d\mu
    \, \leq \, e^{c s} \int_{\R^d} |\nabla \varphi |^2 d\mu, \quad 0 \leq s \leq 1 .
$$
Therefore
$$
 \int_{\R^d} \widehat {Q}_1 \varphi \, d\mu - \int_{\R^d} \varphi \, d\mu
		\, \leq \, \frac {e^{c} -1}{c} \int_{\R^d} |\nabla \varphi |^2 d\mu.
$$
Since
$$
\int_{\R^d} |\nabla \varphi |^2 d\mu \, = \, \int_{\R^d} \varphi (- \mathrm {L} \varphi ) d\mu
		\, = \, \int_{\R^2}  h (- \mathrm{L} )^{-1} h \, d\mu ,
$$
the assertion \eqref {eq.lowerbound1} follows.
\end {proof}

On the basis of Proposition~\ref {prop.lowerbound},
we address the proof of \eqref {eq.lowerboundd=2}.
The first part of the discussion develops in $\R^d$, $d \geq 1$.

The first step is the Kantorovich contraction property under a $CD(0, \infty)$ curvature
condition (cf.~\cite {V09,BGL14}),
which holds in Gauss space for the Mehler kernel $p_t(x,y)$,
\beq \label {eq.contraction}
 \mathrm {W}_2^2 (\mu_n, \mu)   \, \geq \,   \mathrm {W}_2^2 (\mu_n^t, \mu)
\eeq
where we recall that
$d\mu_n^t = f d\mu $, $f(y) = 1 + g(y)$, $ g = g(y) = \frac 1n \sum_{i=1}^n [p_t (X_i,y) - 1]$,
$t >0$.

Next we use the truncation argument on a ball $B_R$ with
radius $R>0$ to be specified later on, and
recall the random variables $X_i^R$, $i = 1, \ldots,n$, with common distribution
$d\mu^R = \frac {1}{\mu(B_R)} \, \mathbbm{1}_{B_R} d\mu$.
Let
$$
{\tilde g} \, = \,  {\tilde g} (y)
	\, = \,  \frac 1n \sum_{i=1}^n \big [ p_t (X_i^R,y) - \E \big ( p_t(X_i^R ,y) \big) \big],
$$
and, for $c>0$,
$$
{\tilde g}_c \, = \,  ({\tilde g}\wedge c)\vee (-c)
   - \int_{\R^d}  \big [({\tilde g} \wedge c)\vee (-c) \big] d\mu
$$
(so that $|{\tilde g}_c| \leq 2c$ and $\int_{\R^d} {\tilde g}_c d\mu = 0$).

In \eqref {eq.lowerbound1} of Proposition~\ref {prop.lowerbound}, choose $h = {\tilde g}_c$.
It holds true that
\beqs \begin {split}
	 \int_{\R^d}  {\tilde g}_c (- \mathrm{L} )^{-1} {\tilde g}_c \, d\mu
   & \, = \, \int_{\R^d}  {\tilde g} (- \mathrm{L} )^{-1} {\tilde g} \, d\mu
      + \int_{\R^d}  ({\tilde g} - {\tilde g}_c) (- \mathrm{L} )^{-1} ({\tilde g} - {\tilde g}_c) d\mu \\
    & \quad \, \,   - 2 \int_{\R^d}  ({\tilde g} - {\tilde g}_c) (- \mathrm{L} )^{-1} {\tilde g} \, d\mu .
\end {split} \eeqs
Therefore, after some algebra, and with $c \leq \frac 12$ for example,
\beq \begin {split} \label {eq.kantorovich2}
 \mathrm {W}_2^2 (\mu^n_t, \mu)
   &\, \geq \, 2 \int_{\R^d} {\tilde g} (- \mathrm{L} )^{-1} g \, d\mu
     -   \frac {e^{2c} -1}{2c}
   \int_{\R^d}  {\tilde g} (- \mathrm{L} )^{-1} {\tilde g} \, d\mu \\
   & \quad \, \, - 2 \int_{\R^d}  ({\tilde g} - {\tilde g}_c) (- \mathrm{L} )^{-1} g \,  d\mu \\
   & \quad \, \, -  2  \int_{\R^d}  ({\tilde g} - {\tilde g}_c)
   (- \mathrm{L} )^{-1} ({\tilde g} - {\tilde g}_c) d\mu
   	 - 4  \, \bigg | \int_{\R^d}  ({\tilde g} - {\tilde g}_c) (- \mathrm{L} )^{-1} {\tilde g} \, d\mu \bigg |. \\
\end {split} \eeq
The three last terms on the right-hand side of \eqref {eq.kantorovich2} are
error terms which may are handled by the exponential
decay \eqref {eq.spectral} in $\mathrm {L}^2(\mu)$.
Indeed, since $ (- \mathrm{L})^{-1}  = \int_0^\infty P_s  ds$,
$$
\int_{\R^d}  ({\tilde g} - {\tilde g}_c) (- \mathrm{L} )^{-1} ({\tilde g} - {\tilde g}_c) d\mu
   \, = \, 2 \int_0^\infty {\big \| P_s ({\tilde g} - {\tilde g}_c) \big\|}_2^2 \, ds
    \, \leq \,  {\|  {\tilde g} - {\tilde g}_c \|}_2^2 .
$$
In the same way,
$$
\int_{\R^d}  ({\tilde g} - {\tilde g}_c) (- \mathrm{L} )^{-1} g \, d\mu
   \, \leq \, {\|  g \|}_2 {\|  {\tilde g} - {\tilde g}_c \|}_2
$$
and
$$
\bigg | \int_{\R^d}  ({\tilde g} - {\tilde g}_c) (- \mathrm{L} )^{-1} {\tilde g} \, d\mu \bigg |
   \, \leq \, {\|  {\tilde g}  \|}_2 {\|  {\tilde g} - {\tilde g}_c \|}_2.
$$
Putting things together, and since
$$
|{\tilde g} - {\tilde g}_c| \, \leq \,  |{\tilde g}| \mathbbm {1}_{\{|{\tilde g}| \geq c\}}
		+ \int_{\R^d}  |{\tilde g}| \mathbbm {1}_{\{|{\tilde g}| \geq c\}} d\mu,
$$
it is deduced from \eqref {eq.contraction} and \eqref {eq.kantorovich2}
that for every $0 < c \leq \frac 12$,
\beq \begin {split} \label {eq.truncation}
  \mathrm {W}_2^2 (\mu^n, \mu)
 	& \, \geq \,  2 \int_{\R^d} {\tilde g} (- \mathrm{L} )^{-1} g \, d\mu
	- \frac {e^{2c} -1}{2c} \int_{\R^d}  {\tilde g} (- \mathrm{L} )^{-1} {\tilde g} \, d\mu \\
    	    & \quad \, \,	- 8 \int_{ \{|{\tilde g}| \geq c\} } |{\tilde g}|^2 d\mu
	    - 8 \big ( {\|  g  \|}_2 + {\|  {\tilde g}  \|}_2 \big)
	    \bigg ( \int_{ \{|{\tilde g}| \geq c\} } |{\tilde g}|^2 d\mu \bigg)^{1/2}. \\
\end {split} \eeq

Next, integrate over the samples
$X_1, \ldots, X_n$ and $X_1^R, \ldots, X_n^R$ the first two terms on the right-hand
side of \eqref {eq.truncation}.
Recalling the definitions of $g$ and ${\tilde g}$, by independence and identical distribution,
$$
\E \bigg ( \int_{\R^d} {\tilde g} (- \mathrm{L} )^{-1} g  \, d\mu \bigg)
  \, = \, \frac 1n \, \int_t^\infty \! \int_{\R^d}
      \E \Big ( \big [ p_t (X_1^R,y) - \E \big ( p_t(X_1^R ,y) \big) \big] p_s(X_1,y ) \Big) d\mu(y) ds.
$$
By definition of $X_1^R$,
\beqs \begin {split}
 \E \Big ( \big [ p_t (X_1^R,y) - \E & \big ( p_t(X_1^R ,y) \big)  \big] p_s(X_1,y ) \Big) \\
 & \, = \, \E \Big ( \mathbbm{1}_{\{X_1 \in B_R\}}
 \big [ p_t (X_1,y) - \E \big ( p_t(X_1^R ,y) \big) \big] p_s(X_1,y ) \Big) \\
 & \quad \, \, + \E \Big (\mathbbm{1}_{\{X_1 \notin B_R\}}
 \big [ p_t (Z_1,y) - \E \big ( p_t(X_1^R ,y) \big) \big] p_s(X_1,y ) \Big) \\
 & \, = \, \E \Big ( \mathbbm{1}_{\{X_1 \in B_R\}}
 \big [ p_t (X_1,y) - \E \big ( p_t(X_1^R ,y) \big) \big] p_s(X_1,y ) \Big) \\
\end {split} \eeqs
since $Z_1$ is independent of $X_1$ and with the same law as $X_1^R$.
Hence, after integration in $d\mu(y)$ and the semigroup property,
\beqs \begin {split}
\E \bigg ( \int_{\R^d} {\tilde g} & (- \mathrm{L} )^{-1} g \, d\mu \bigg) \\
   & \, = \,  \frac {\mu(B_R)}{n} \, \int_t^\infty  \bigg [ \int_{\R^d} p_{t+s} (x,x) d\mu^R(x)
         - \int_{\R^d} \!\int_{\R^d} p_{t+s} (x,x') d\mu^R(x) d\mu^R(x') \bigg] ds \\
           & \, = \,  \frac {\mu(B_R)}{n} \, \int_{2t}^\infty  \bigg [ \int_{\R^d} p_s (x,x) d\mu^R(x)
         - \int_{\R^d} \!\int_{\R^d} p_s (x,x') d\mu^R(x) d\mu^R(x') \bigg] ds \\
\end {split} \eeqs
In the same way,
\beqs \begin {split}
\E \bigg ( \int_{\R^d} {\tilde g} & (- \mathrm{L} )^{-1} {\tilde g} \, d\mu \bigg) \\
   & \, = \,  \frac 1n \, \int_t^\infty \! \int_{\R^d}
      \E \Big ( \big [ p_t (X_1^R,y) - \E \big ( p_t(X_1^R ,y) \big) \big] p_s(X_1^R,y ) \Big) d\mu(y) ds \\
   & \, = \,  \frac 1n \, \int_{2t}^\infty  \bigg [ \int_{\R^d} p_s (x,x) d\mu^R(x)
         - \int_{\R^d} \!\int_{\R^d} p_s (x,x') d\mu^R(x) d\mu^R(x') \bigg] ds.
\end {split} \eeqs

As a consequence, if $c>0$ is small enough and $\mu(B_R)$ close to $1$,
\beqs \begin {split}
\E \bigg ( 2 \int_{\R^d} {\tilde g}&  (- \mathrm{L} )^{-1} g \, d\mu
	- \frac {e^{2c} -1}{2c} \int_{\R^d}  {\tilde g} (- \mathrm{L} )^{-1} {\tilde g} \, d\mu  \bigg) \\
	& \, \geq \, \frac {1}{2n}
	\bigg (  \int_{2t}^\infty  \bigg [ \int_{\R^d} p_s (x,x) d\mu^R(x)
         - \int_{\R^d} \!\int_{\R^d} p_s (x,x') d\mu^R(x) d\mu^R(x') \bigg] ds \bigg).
\end {split} \eeqs
Also, by the spectral gap inequality \eqref {eq.spectral},
\beqs \begin {split}
\int_{\R^d} \!\int_{\R^d} p_s (x,x') d\mu^R(x) d\mu^R(x')
		&\, = \, \frac {1}{\mu(B_R)^2} \int_{\R^d}
			\mathbbm {1}_{B_R} P_s \mathbbm {1}_{B_R} \, d\mu  \\
		&\, = \, 1 + \frac {1}{\mu(B_R)^2} \int_{\R^d}
			\mathbbm {1}_{B_R} P_s \big (\mathbbm {1}_{B_R} - \mu(B_R) \big) d\mu  \\
		&\, \leq \, 1 + \frac {1 - \mu(B_R)}{\mu(B_R)} \, e^{-s} \\
		&\, \leq \, 1 + 2\, e^{-s}
\end {split} \eeqs
provided that $\mu(B_R) \geq \frac 12$. As a conclusion at this stage,
\beq \begin {split} \label {eq.expectation}
 \E \big ( \mathrm {W}_2^2 (\mu^n, \mu)  \big )
	& \, \geq \,
	\, \frac {1}{2n} \int_{2t}^\infty  \! \int_{\R^d} \big [ p_s (x,x) - 1 \big] d\mu^R(x) ds
         - \frac 1n \\
             & \quad \, \,  - 8 \, \E \bigg ( \int_{ \{|{\tilde g}| \geq c\} } |{\tilde g}|^2 d\mu \bigg)
            - 8 \,  \E \bigg ( \big ( {\|  g  \|}_2 + {\|  {\tilde g}  \|}_2 \big)
	   \bigg ( \int_{ \{|{\tilde g}| \geq c\} } |{\tilde g}|^2 d\mu \bigg)^{1/2} \, \bigg)  . \\
\end {split} \eeq

\medskip

The final part of the proof will be to take care of the correction terms
on the right-hand side of the preceding \eqref {eq.expectation}. To this task,
we develop some (crude) bounds on the Mehler kernel $p_t(x,y)$ of Section~\ref {sec.2}.
Consider for each $y \in \R^d$, $t>0$ and $q \geq 1$,
$$
\int_{B_R} p_t (x,y)^q d\mu(x).
$$
After translation and a change of variable,
$$
\int_{B_R} p_t  (x,y)^q  d\mu(x) \\
	 \, = \, \frac {1}{(1 -a^2)^{qd/2}} \, e^{\tau (q-1)a\frac {|y|^2}{2}}
	    \int_{B(-\tau y, R)}
		\, e^{- \frac {qa}{\tau (1-a^2)} \,\frac {|x|^2}{2}} \frac {dx}{(2\pi)^{d/2}}
$$
where $a = e^{-t}$ and $\tau = \frac {q a}{1+(q -1)a^2}$.

Note that $ \tau (q-1)a \leq q$,
and that $\tau \geq \frac 12$ at least provided that $a $ is close to one which we may assume.
Then, if $|y| \geq 4R$ and $ x \in B(-\tau y, R)$, we have $|x| \geq \frac {|y|}{4}$. Hence,
whenever $|y| \geq 4R$,
$$
\int_{B_R} p_t (x,y)^q d\mu(x) \, \leq \, \frac {1}{(1 -a^2)^{(q -1)d/2}} \,
		e^{- \big( \frac {1}{32(1-a^2)} - \frac q2\big) |y|^2} .
$$
Otherwise, that is when $|y| \leq 4R$,
$$
\int_{B_R} p_t (x,y)^q d\mu(x) \, \leq \,
	\frac {1}{(1 -a^2)^{(q -1)d/2}} \, e^{\frac q2 |y|^2} .
$$

Recall
$$
{\tilde g} \, = \,  {\tilde g} (y)
	\, = \,  \frac 1n \sum_{i=1}^n \big [ p_t (X_i^R,y) - \E \big ( p_t(X_i^R ,y) \big) \big] .
$$
By Rosenthal's inequality \eqref {eq.rosenthal},
for any $q \geq 2$ there exists $C_q >0$ only depending on $q$ such that
\beqs \begin {split}
 \E \big ( \big | {\tilde g}(y) \big |^q \big)
     & \, \leq \, C_q \bigg ( \frac {1}{n^{q-1}} \, \E \big ( p_t (X_1^R,y)^q \big)
             + \frac {1}{n^{q/2}} \,  \Big [ \E \big ( p_t (X_1^R,y)^2\big) \Big ]^{q/2} \bigg) \\
      & \, \leq \, 2^{\frac q2 } C_q \bigg ( \frac {1}{n^{q-1}} \, \int_{B_R} p_t (x,y)^q d\mu(x)
             + \frac {1}{n^{q/2}} \,  \bigg [ \int_{B_R} p_t (x,y)^2 d\mu(x) \bigg ]^{q/2} \bigg) \\
\end {split} \eeqs
where it is assumed that $\mu (B_R) \geq \frac 12$.

In the following $q \geq 2$ is fixed. Then $t>0$ may be chosen small enough (in terms
of $q$ but independently of $n$) such that $\tau  \geq \frac 12$ and
$ \frac {1}{32(1-a^2)} - \frac q2 \geq 0$ (for example). By the previous step,
$$
\int_{\R^d} \bigg (\int_{B_R} p_t (x,y)^q d\mu(x) \bigg) d\mu(y)
	\, \leq \, \frac {1}{ (1 -a^2)^{(q -1)d/2}} \, \big ( 1 +  e^{8qR^2} \big) .
$$
In the same way,
$$
\int_{\R^d} \bigg (\int_{B_R} p_t (x,y)^2 d\mu(x) \bigg)^{q/2} d\mu(y)
	\, \leq \, \frac {1}{ (1 -a^2)^{qd/4}} \, \big ( 1 +  e^{4qR^2} \big) .
$$
For simplicity (in order not to carry the two preceding
expressions with $q-1$ and $\frac q2$),
assume in the following that $(1-a^2)^{d/2} n \geq 1$.
Therefore, using that $q-1 \geq \frac q2$,
\beq \label {eq.rosenthal2}
\int_{\R^2}  \E \big ( \big | {\tilde g}(y) \big |^q \big) d\mu (y)
		\, \leq \,  \frac {2^qC_q}{ [(1 -a^2)^{d/2}n]^{q/2}} \, \big ( 1 +  e^{8qR^2} \big) .
\eeq

We use the preceding bounds to control the error term
\beq \label {eq.error}
 \mathrm {Er} \, = \,  \E \bigg ( \int_{ \{|{\tilde g}| \geq c\} } |{\tilde g}|^2 d\mu \bigg)
            +   \E \bigg ( \big ( {\|  g  \|}_2 + {\|  {\tilde g}  \|}_2 \big)
	   \bigg ( \int_{ \{|{\tilde g}| \geq c\} } |{\tilde g}|^2 d\mu \bigg)^{1/2} \, \bigg)
\eeq
of \eqref {eq.expectation}.
By repeated use of the Young and H\"older inequalities, the latter is bounded from above
for any $\delta > 0$ and $\alpha >1$ by
$$
\delta \, \Big [ \E \big (  {\|  g \|}_2^2 \big) +  \E \big (  {\|  {\tilde g}  \|}_2^2 \big) \Big]
 + \frac {1 + 2 \delta}{2 c^{2(\alpha -1) }\delta} \,
 		 \int_{\R^d} \E \big (  \big | {\tilde g} (y)  \big |^{2\alpha} \big )d\mu(y) .
$$
Since $p_t(x,x) = \frac {1}{1-a^2} \,e^{\frac {a}{1+a} |x|^2}$, again with $\mu (B_R) \geq \frac 12$,
\beqs \begin {split}
\E \big (  {\|  {\tilde g}  \|}_2^2 \big)
	& \, = \, \frac 1n \int_{\R^d} \Big [ \E \big( p_t(X_1^R,y)^2 \big)
			- \E \big( p_t(X_1^R,y) \big)^2 \Big] d\mu(y) \\
	& \, \leq \, \frac 1n \, \int_{\R^d} p_t (x,x) d\mu^R (x)  \\
	&\, \leq \,  \frac {1}{n \mu(B_R)(1 - a)^d} \\
	&	\, \leq \,  \frac {2}{n (1 - a)^d} \, .
\end {split} \eeqs
Similarly
$$
\E \big (  {\|  g  \|}_2^2 \big)
		\, \leq \, \frac 1n \, \int_{\R^2} p_t (x,x) d\mu (x)
		\, \leq \, \frac {1}{n (1 - a)^d} \, .
$$
On the other hand, \eqref {eq.rosenthal2} with $q = 2 \alpha $ yields
$$
\int_{\R^d} \E \big (  \big | {\tilde g} (y)  \big |^{2\alpha} \big) d\mu(y)
	\, \leq \, \frac {4^\alpha C_{2\alpha}}{ [(1 -a^2)^{d/2}n]^{\alpha}} \, \big ( 1 +  e^{16\alpha R^2} \big) .
$$
Hence, for any $0 < \delta \leq 1$ and $\alpha >1$, the error term \eqref {eq.error} satisfies
$$
	\mathrm {Er}  \, \leq \, \frac {3\delta}{n (1 - a)^d} +
		\frac {4^{\alpha +1} C_{2\alpha}}{  c^{2(\alpha -1)} \delta [(1 -a^2)^{d/2}n]^{\alpha}}
		  \big ( 1 +  e^{16\alpha R^2} \big) .
$$
Therefore, from \eqref {eq.expectation},
\beq \begin {split} \label {eq.expectationerror}
 \E \big ( \mathrm {W}_2^2 (\mu^n, \mu)  \big )
	& \, \geq \,
	\, \frac {1}{2n} \int_{2t}^\infty  \! \int_{\R^d} \big [ p_s (x,x) - 1 \big] d\mu^R(x) ds
         - \frac 1n \\
             & \quad \, \,  - \frac {24\delta}{n (1 - a)^d} -
		\frac {4^{\alpha +3} C_{2\alpha}}{  c^{2(\alpha -1)} \delta [(1 -a^2)^{d/2}n]^{\alpha}}
		  \big ( 1 +  e^{16\alpha R^2} \big)   . \\
\end {split} \eeq

\bigskip

In this last step, we fix the various parameters involved in the previous analysis.
Basically, $ t \approx \frac {1}{n^\varepsilon}$ and $R \approx \varepsilon \sqrt {\log n}$
for some small $\varepsilon >0$, and $\alpha >1$ is chosen large enough.
Take for example $ t = \frac {1}{n^{1/d}}$ and
$R^2 = \frac {1}{64} \log n$. Then, for $n$ large enough, the necessary conditions
on $a = e^{-t}$ or $\mu(B_R)$ are fulfilled. After some details, the choice of $\delta = \frac 1n$
and $\alpha = 8$ in \eqref {eq.expectationerror} yields that
$$
 \E \big ( \mathrm {W}_2^2 (\mu^n, \mu)  \big )
	 \, \geq \, \frac {1}{2n} \int_{2t}^\infty  \! \int_{\R^d} \big [ p_s (x,x) - 1 \big] d\mu^R(x) ds
		- O \Big (\frac 1n \Big).
$$

From the analysis of the upper bound in \cite{L17}, it is known that
as $ t < \! <  \frac {1}{R^2}$, for some $\rho >0$,
$$
\int_{2t}^\infty  \! \int_{\R^2} \big [ p_s (x,x) - 1 \big] d\mu^R(x) ds
   \, \geq \, \rho R^2 \log \Big( \frac 1t \Big)
$$
when $d=2$, and also
$$
\int_{2t}^\infty  \! \int_{\R} \big [ p_s (x,x) - 1 \big] d\mu^R(x) ds
   \, \geq \, \rho \log (R^2)
$$
when $d=1$. Therefore, for the preceding choices of $t$ and $R$, this establishes the claim
\eqref {eq.lowerboundd=2}, as well as the  $\frac {\log \log n}{n}$ lower bound in \eqref {eq.gaussian1}.
The announced proof is complete.

\vskip 8mm

\font\tenrm =cmr10  {\tenrm

\parskip 0mm

\noindent Institut de Math\'ematiques de Toulouse, Universit\'e de Toulouse -- Paul-Sabatier, F-31062 Toulouse, France \&  Institut Universitaire de France, ledoux@math.univ-toulouse.fr

\medskip

\noindent School of Mathematical Sciences, Fudan University, Shanghai 200433, People's Republic of China,15110840006@fudan.edu.cn \& Institut de Math\'ematiques de Toulouse,
Universit\'e de Toulouse -- Paul-Sabatier, F-31062 Toulouse, France,
zhu@math.univ-toulouse.fr

}


\vskip 7mm





\begin{thebibliography}{9}


\bibitem {AKT84}
M. Ajtai, J. Koml\'os, G. Tusn\'ady. On optimal matchings. \textit {Combinatorica}~4, 259--264 (1984).

\bibitem {AST19}
L. Ambrosio, F. Stra, D. Trevisan. A PDE approach to a 2-dimensional matching problem (2016).
\textit {Probab. Theory Related Fields}~173, 433-478 (2019).

\bibitem {B87}
D. Bakry. \'Etude des transformations de Riesz sur les variétés riemanniennes à courbure de Ricci minorée.
\textit {Séminaire de Probabilités XXI, Lecture Notes in Math.}~1247, 137--172. Springer (1987).

\bibitem {BGL14}
D. Bakry, I. Gentil, M. Ledoux. \textit {Analysis and geometry of Markov
diffusion operators.} Grundlehren der mathematischen Wissenschaften~348.
Springer (2014).

\bibitem {BL16}
S. Bobkov, M. Ledoux. One-dimensional empirical measures, order statistics, and Kantorovich
transport distances (2016). To appear in \textit {Memoirs Amer. Math. Soc.}

\bibitem {DSS13}
S. Dereich, M. Scheutzow, R. Schottstedt. Constructive quantization: approximation by empirical
measures. \textit {Ann. Inst. Henri Poincaré Probab. Stat.}~49, 1183--1203 (2013).

\bibitem {FG15}
N. Fournier, A. Guillin. On the rate of convergence in Wasserstein distance
of the empirical measure. \textit {Probab. Theory Related Fields}~162, 707--738 (2015).

\bibitem {L17}
M. Ledoux. On optimal matching of Gaussian samples. \textit {Zap. Nauchn. Sem. S.-Peterburg. Otdel. Mat. Inst. Steklov. (POMI) 457, Veroyatnost' i Statistika.} 25, 226--264 (2017).

\bibitem {L18}
M. Ledoux. On optimal matching of Gaussian samples II (2018).

\bibitem {M84}
P.-A. Meyer. Transformations de Riesz pour les lois gaussiennes
\textit {Séminaire de Probabilités XV, Lecture Notes in Math.}~1059, 179--193. Springer (1984).

\bibitem {P18}
R. Peyre. Comparison between $\mathrm{W}_2$ distance and
$\dot{\mathrm {H}}^{-1}$ norm, and localization
of Wasserstein distance. \textit {ESAIM Control Optim. Calc. Var.}~24, 1489–1501 (2018).

\bibitem {R70}
H. P. Rosenthal. On the subspaces of $L^p$ $(p > 2)$ spanned by sequences of independent
random variables. \textit {Israel J. Math.}~8, 273--303 (1970).

\bibitem {S15}
F. Santambrogio. \textit {Optimal Transport for Applied Mathematicians.}
Progress in Nonlinear Differential Equations and Their Applications. Birkh\"auser (2015).

\bibitem {T14}
M. Talagrand.
\textit {Upper and lower bounds of stochastic processes.}
Modern Surveys in Mathematics~60. Springer-Verlag (2014).

\bibitem {T18}
M. Talagrand. Scaling and non-standard matching theorems. \textit {Comptes Rendus
Acad. Sciences Paris, Mathématique}~356, 692--695 (2018).

\bibitem {V09}
C. Villani. \textit {Optimal transport. Old and new.}
Grundlehren der mathematischen Wissenschaften~338. Springer (2009).

\bibitem {Y92}
J. Yukich. Some generalizations of the Euclidean two-sample matching problem.
\textit {Probability in Banach spaces}~8, Progr. Probab.~30, 55--66. Birkh\"auser (1992).


\end{thebibliography}
\end{document}